\documentclass{article} 
\usepackage{microtype}

\usepackage{amsmath}\allowdisplaybreaks
\usepackage{amsfonts,bm}
\usepackage{amssymb}

\def\eqref#1{equation~(\ref{#1})}

\def\1{\bf{1}}

\newcommand{\Norm}[1]{\left\| #1 \right\|}
\newcommand{\norm}[1]{\left\| #1 \right\|}

\def\vzero{{\bf{0}}}
\def\vone{{\bf{1}}}
\def\vxi{{\bm{\xi}}}

\def\va{{\bf{a}}}

\def\vc{{\bf{c}}}

\def\vg{{\bf{g}}}

\def\vu{{\bf{u}}}
\def\vv{{\bf{v}}}
\def\vw{{\bf{w}}}
\def\vx{{\bf{x}}}
\def\vy{{\bf{y}}}
\def\vz{{\bf{z}}}

\def\fD{{\mathcal{D}}}

\def\fO{{\mathcal{O}}}

\def\BE{{\mathbb{E}}}
\def\BB{{\mathbb{B}}}

\def\BR{{\mathbb{R}}}
\def\BS{{\mathbb{S}}}

\def\mG {{\bf G}}

\def\mI {{\bf I}}

\def\mP {{\bf P}}

\def\mV {{\bf V}}
\def\mW {{\bf W}}
\def\mX {{\bf X}}

\def\mZ {{\bf Z}}

\def\FG{{\rm FastGossip}}

\usepackage{amsthm}
\usepackage{etoolbox}

\makeatletter
\def\Ddots{\mathinner{\mkern1mu\raise\p@
\vbox{\kern7\p@\hbox{.}}\mkern2mu
\raise4\p@\hbox{.}\mkern2mu\raise7\p@\hbox{.}\mkern1mu}}
\makeatother

\makeatletter
\newcommand*{\rom}[1]{\expandafter\@slowromancap\romannumeral #1@}
\makeatother

\def\DOCS {{DOC$^2$S}}

\def\vnabla {{{\bm \nabla}}}
\def\vxi {{{\bm \xi}}}

\newcommand{\fnormsq}[1]{\left\| #1 \right\|_F^2}
\newcommand{\fnorm}[1]{\left\| #1 \right\|_F}

\usepackage{subcaption}

\usepackage{hyperref}
\usepackage{url}
\usepackage{mathrsfs}
\usepackage[letterpaper,margin=1.0in]{geometry}
\usepackage[utf8]{inputenc} 
\usepackage[T1]{fontenc}    
\usepackage{booktabs}       
\usepackage{amsfonts}       
\usepackage{nicefrac}       
\usepackage{microtype}      
\usepackage{xcolor}
\usepackage[numbers,sort,compress]{natbib}
\usepackage{graphicx}
\usepackage{microtype}
\usepackage{enumitem}
\usepackage{algorithm}
\usepackage{algorithmic}
\usepackage[capitalize,noabbrev]{cleveref}
\usepackage[textsize=tiny]{todonotes}

\theoremstyle{plain}
\newtheorem{theorem}{Theorem}
\newtheorem{proposition}{Proposition}
\newtheorem{lemma}{Lemma}
\newtheorem{corollary}[theorem]{Corollary}
\theoremstyle{definition}
\newtheorem{definition}{Definition}
\newtheorem{assumption}{Assumption}
\theoremstyle{remark}
\newtheorem{remark}{Remark}
\newcommand{\gradnormdelta}[2]{\norm {\nabla #1 (#2) }_\delta }

\newcommand{\un}{{\rm Unif}}

\newcommand{\finitesum}[2]{\sum_{#1=1}^#2}

\newcommand{\vclip}{{\bf clip}}
\newcommand{\mClip}{{\bf Clip}}

\title{\vskip-0.5cm Decentralized Nonsmooth Nonconvex Optimization with \\ Client-Data Sampling}

\date{}
\author{\vspace{-4cm} Cheng Chen \quad\quad\quad Xinyan Chen \quad\quad\quad Weiguo Gao \quad\quad\quad Luo Luo}

\begin{document}

\maketitle

\begin{abstract}
This paper studies decentralized nonsmooth nonconvex optimization with Lipschitz-continuous local functions. 
We propose an efficient stochastic first-order method with a novel client-data sampling strategy in which the mini-batch sizes of the stochastic gradients on all clients follow a multinomial distribution.
We show that our method can find a $(\delta,\epsilon)$-Goldstein stationary point on each client with the best-known dependence on 
$\epsilon$~and $\delta$ in both computational and communication complexity.
Moreover, our results tighten the upper complexity bounds with respect to the Lipschitz parameter, the variance upper bound,  the number of clients, and the spectral gap.
We further conduct numerical experiments to demonstrate the empirical advantages of the proposed method.

\end{abstract}

\section{Introduction}
The large scale nonsmooth nonconvex optimization  covers many applications in fields such as machine learning \citep{nair2010rectified,xiao2023Nonsmooth}, statistics \citep{fan2001variable,zhang2010nearly,zhang2010analysis}, and economics \citep{duffie2010dynamic,stadtler2014supply}.
In this paper, we focus on the decentralized stochastic optimization problem 
\begin{align}\label{eq:global}
   \min_{\vx \in \mathbb{R}^d} f(\vx) \triangleq \frac{1}{n} \sum_{i=1}^{n} f_i(\vx)  
\end{align}
over the network with $n$ clients, where the local function at the $i$th client has the form of
\begin{align}\label{eq:local-func}
f_i(\vx) \triangleq \mathbb{E}_{\vxi_i \sim \mathcal{D}_i}[F_i(\vx; \vxi_i)]    
\end{align}
such that the $f_i:\BR^d\to\BR$ is Lipschitz continuous but possibly nonconvex nonsmooth and the random index $\vxi_i\in\Xi_i$ follows the distribution $\mathcal{D}_i$. 
It is well known that achieving approximate stationary points in terms of the classical Clarke subdifferential \citep{clarke1990optimization} for the general Lipschitz continuous function is intractable \citep{jordan2022complexity,kornowski2022oracle,tian2022no,zhang2020complexity}. 
Instead, we typically target to find $(\delta,\epsilon)$-Goldstein stationary points \citep{zhang2020complexity}. 
This criterion suggests studying the convex hull of Clarke subdifferential at points in the $\delta$-radius neighborhood of the given point. 

The stochastic optimization methods for finding $(\delta,\epsilon)$-Goldstein stationary points in non-distributed setting have been widely studied in recent years \citep{chen2023faster,cutkosky2023optimal,davis2022gradient,kornowski2023algorithm,lin2022gradient,tian2022finitetime,zhang2020complexity}.
Specifically, \citet{tian2022finitetime,zhang2020complexity} proposed the (perturbed) stochastic interpolated normalized gradient descent with the stochastic first-order oracle (SFO) complexity of~$\fO(\delta^{-1}\epsilon^{-4})$.
In a seminal work, \citet{cutkosky2023optimal} established the conversion from nonsmooth nonconvex optimization to online learning, achieving the SFO complexity of $\fO(\delta^{-1}\epsilon^{-3})$.
They also extends the lower bound of \citet{arjevani2023lower} to show their SFO complexity is optimal.

In decentralized setting, we have to consider the consensus error among the variables on different clients in the network.
The popular technique of gradient tracking can successfully bound the consensus error in for smooth optimization problems  \citep{nedic2009distributed,qu2017harnessing,shi2015extra}, while it cannot be directly used in the nonsmooth setting since the Lipschitz continuity of the gradient (subgradient) may not hold.  
\citet{kovalev2024lower,lan2020communication} proposed efficient decentralized algorithms based on the primal-dual framework for the nonsmooth objective but only limit to the convex problem.
A natural idea for decentralized nonsmooth optimization is using the randomized smoothing to establish the smooth surrogate for the original problem, which works for both the convex \citep{scaman2018optimal} and the nonconvex settings \citep{lin2023decentralized,sahinoglu2024online}.
Notably, \citet{sahinoglu2024online} proposed multi-epoch decentralized online learning (ME-DOL) method for decentralized stochastic stochastic nonsmooth nonconvex optimization, which incorporates the decentralized online mirror descent \citep{shahrampour2017distributed} into the online-to-nonconvex conversion \citep{cutkosky2023optimal,kornowski2023algorithm}.
The theoretical analysis shows that ME-DOL can find $(\delta,\epsilon)$-Goldstein stationary points under the assumption of bounded second moment for the stochastic gradient \citep{sahinoglu2024online}.

\begin{table*}[t]
{\centering
\caption{
We present the upper complexity bounds of our methods and related work for finding $(\delta,\epsilon)$-Goldstein stationary points in stochastic decentralized nonsmooth nonconvex optimization problem.
Here, the sample complexity refers to the total local stochastic first-order oracle calls on all of the $n$ clients.
}
\label{sample-table1}
\resizebox{\linewidth}{!}{\begin{tabular}{cccc}
\hline
 Methods & $^\dag$ME-DOL \citep{sahinoglu2024online} & \DOCS~(Theorem 1) \\
 \hline \addlinespace
 Sample Complexity & $\fO\bigg(
 \dfrac{L+\Delta}{\gamma^{2}\delta\epsilon^{3}}\big(
 G(L+G)n^2+(G^2+L^2d)\gamma n\sqrt{n}\,\big)\bigg)$ & $\fO\bigg(\dfrac{(bL^2+\sigma^2)\Delta}{\delta\epsilon^3}\bigg)$ \\ \addlinespace
 Computation Rounds & $\fO\bigg(
 \dfrac{L+\Delta}{\gamma^{2}\delta\epsilon^{3}}\big(
 G(L+G)n+(G^2+L^2d)\gamma\sqrt{n}\,\big)\bigg)$ & $\tilde\fO\bigg(\dfrac{(bL^2+\sigma^2)\Delta}{n\delta\epsilon^3}\bigg)$ \\ \addlinespace
 Communication Rounds & $\fO\bigg(
 \dfrac{L+\Delta}{\gamma^{2}\delta\epsilon^{3}}\big(
 G(L+G)n+(G^2+L^2d)\gamma\sqrt{n}\,\big)\bigg)$ & $\tilde\fO\bigg(\dfrac{(L^2+b^{-1}\sigma^2)\Delta}{\sqrt{\gamma}\delta\epsilon^{3}}\bigg)$ \\ \addlinespace
 \hline
\end{tabular}}}
\vskip0.1cm
\begin{itemize}[leftmargin=0.2cm]
\item[$^\dag$\!\!]{\footnotesize ME-DOL \citep[Assumption 4]{sahinoglu2024online} additionally assumes the stochastic gradient has bounded second moment, i.e., there exists constant $G>0$ such that~$\BE[\|\nabla F_i(\vx,\vxi_i)\|^2]\le G^2$ for all $\vx\in\BR^d$.
}
\end{itemize} \vspace{-1cm}
\end{table*}

It is worth noting that the objective in distributed optimization problem (\ref{eq:global}) naturally has the finite-sum structure in the view of local functions $\{f_i\}_{i=1}^n$.
This motivates us to design the partial participated methods, which performs the client sampling during iterations and only executes the computation on the selected clients \citep{chen2020optimal,maranjyan2022gradskip,mishchenko2022proxskip}.
Some recent works \citep{bai2024complexity,liudecentralized,luo2022optimal,chen2025near} studied  partial participated methods by considering the balance among the first-order oracle complexity, the computation rounds, and the communication rounds in decentralized optimization. 
However, these results heavily depend on the smoothness assumptions.
To the best of our knowledge, all existing methods for decentralized nonsmooth optimization \citep{chen2020optimal,kovalev2024lower,lan2020communication,lin2023decentralized,sahinoglu2024online,wang2023decentralized,zhang2024decentralized} require all clients accessing their local oracle in per computation rounds,
which limits the sampling efficiency.
 
In this paper, we propose \underline{d}ecentralized \underline{o}nline-to-nonconvex \underline{c}onversion with \underline{c}lient-data \underline{s}ampling (\DOCS),  incorporating a novel client-data sampling scheme and Chebyshev acceleration~\citep{ye2023multi,arioli2014chebyshev}.
The proposed client-data sampling scheme repeats the sampling procedure $b$ times, where each of them first samples a client and then draws a stochastic gradient from the selected client. This two-stage sampling strategy enables \DOCS~to perform mini-batch update with unfixed batch sizes on different clients,
resulting the sharper dependence on the Lipschitz parameter and variance upper bound. 
On the other hand, the Chebyshev acceleration allows \DOCS~significantly reduce the consensus error in each round and therefore improve the dependence on the number of clients $n$ and the dimensionality $d$ in the upper complexity bounds of ME-DOL \citep{sahinoglu2024online}. 
We show that \DOCS~can achieve an $(\delta,\epsilon)$-Goldstein stationary point on each client with the total sample complexity of $\fO((bL^2+\sigma^2)\delta^{-1}\epsilon^{-3}\Delta)$, 
the computation rounds of $\fO((L^2+b^{-1}\sigma^2)\delta^{-1}\epsilon^{-3}\Delta)$,
and the communication rounds of $\fO((bL^2+\sigma^2)\gamma^{-1/2}\delta^{-1}\epsilon^{-3}\Delta)$, 
where~$b$ is the overall batch size in each iteration, $L$ is the Lipschitz constant of local functions, $\sigma^2$ is the upper bound of the stochastic gradient variance, $\gamma$ is the spectral gap of the mixing matrix, and~$\Delta$ is the initial function value gap.
These upper complexity bounds indicate the mini-batch size $b$ in our algorithm provides a tradeoff between the computation complexity and the communication complexity.
Specifically, even by taking $b=1$, all of our upper bounds are tighter than those of ME-DOL \citep{sahinoglu2024online}.
We present the comparison of theoretical results of our methods and related work in Table~\ref{sample-table1}.

\section{Preliminaries}
In this section, we formalize our problem setting and introduce the background of nonsmooth analysis.

\subsection{Notation and Assumptions}

We use the bold lowercase letter with the subscript to present the local variable, e.g., the row vector~$\vx_i\in \mathbb{R}^{1\times d}$ presents the local variable on the $i$th client. 
The bold lowercase letter with a bar presents the average of local
variables, e.g., $\bar\vx=\frac{1}{n}\sum_{i=1}^n\vx_i$.
We also use bold uppercase letters to denote matrices that consists of local variables, e.g.,
$\mX=\begin{bmatrix}
    \vx_{1} ; \dots ; \vx_{n}
\end{bmatrix}    
\in \mathbb{R}^{n \times d}$.
For the ease of presentation, we allow the input of the function and its gradient to be organized as either a column vector or a row vector. 
We denote $\|\cdot\|$ and $\|\cdot\|_2$ as the Frobenius norm and the spectral norm of the matrix, respectively, also the Euclidean norm of the vector.
We let $\vone_n=[1,\dots,1]^\top\in\BR^n$ and $\mI$ be the identity matrix.
The notation $\mathrm{conv}(\cdot)$ denotes the convex hull of given set.
Additionally, we use notations $\BB^{d}(\vx,\delta)$ and $\BS^{d-1}$ to present the Euclidean ball centered at $\vx\in\BR^d$ with radius~$\delta>0$ and the unit sphere centered at the origin, respectively. 

We impose following assumptions for formulations (\ref{eq:global}) and (\ref{eq:local-func}).
\begin{assumption}\label{assumption:Lipschitz}
We suppose each local function $f_i:\BR^d\to\BR$ is $L$-Lipschitz continuous, i.e., there exists $L>0$ such that 
$|f_i(\vx)-f_i(\vy)|\leq L||\vx-\vy||$ for all $\vx,\vy\in\BR^d$ and $i\in[n]$.
\end{assumption}

\begin{assumption}
We suppose the objective function $f:\BR^d\to\BR$ is lower bounded by $f^*$, i.e., it holds~$f^* \triangleq \inf_{\vx\in\BR^d} f(\vx) > -\infty$. 
\label{assumption:suboptimality}
\end{assumption}
Note that Rademacher's theorem \citep{evans2018measure} states  a Lipschitz continuous function is differentiable almost everywhere. Thus, the LSFO is well-defined almost everywhere under Assumption~\ref{assumption:Lipschitz}.

We make the following assumption for the local stochastic first-order oracle (LSFO).

\begin{assumption}\label{assumption:gradientoracle}
We suppose the algorithms can access the local function $f_i:\BR^d\to\BR$ via the LSFO consisting of the local gradient estimator $F_i:\BR^d\times\Xi_i\to\BR$ and the random index $\vxi_i\sim\fD_i$ such that~$\mathbb{E}_{\vxi_i\sim\fD_i}[\nabla F_i(\vx;\vxi_i)] = \nabla f_i(\vx)$
and 
$\mathbb{E}_{\vxi_i\sim\fD_i}[\norm{ \nabla F_i(\vx;\vxi_i) - \nabla f_i(\vx) }^2 ] \le \sigma^2$ 
for some $\sigma\geq 0$.
\end{assumption}

We aim for all $n$ clients in the network to collaborate in solving stochastic decentralized optimization problem~(\ref{eq:global}).
We use the doubly stochastic matrix $\mP=[p_{ij}]\in\BR^{n\times n}$ to describe the topology of the network.
Specifically, the communication step at the $i$th client is built upon the weighted average~$\vx_i^+=\sum_{j=1}^n p_{ij}\vx_j$, where  $\vx_j$ is the local variable on the $j$th client.
We impose the following standard assumption in decentralized optimization for the matrix $\mP\in\BR^{n\times n}$ \citep{schmidt2017minimizing,scaman2018optimal}.

\begin{assumption}\label{assumption:P}
We suppose that the mixing matrix $\mP\in \mathbb{R}^{n \times n}$  associated with the network satisfies:
(a) The matrix $\mP\in\BR^{n\times n}$ is symmetric and holds $p_{ij} \ge 0$ for all $i,j\in[n]$;~
(b) It holds $p_{ij} \neq 0$ if and only if clients $i$ and $j$ are connected or $i = j$;~
(c) It holds $\vzero \preceq \mP \preceq \mI$ and $\mP^\top \vone_n= \mP\vone_n = \vone_n$.
\end{assumption}

Under Assumption \ref{assumption:P}, the largest eigenvalue of the mixing matrix $\mP\in\BR^{n\times n}$ is one.
Consequently, we define the spectral gap of $\mP\in\BR^{n\times n}$ as
$\gamma = 1-\lambda_2(\mP) \in (0,1]$,
where $\lambda_2(\mP)$ is the second largest eigenvalue of $\mP$.

\subsection{Goldstein Stationary Points}

We present the notion of Clarke subdifferential \citep{clarke2008nonsmooth} and its relaxation Goldstein subdifferential \citep{goldstein1977optimization}
for the Lipschitz continuous objective in the nonconvex nonsmooth problem as follows.

\begin{definition}[\citet{clarke2008nonsmooth}]\label{dfn:Clk}
The Clarke subdifferential of a Lipschitz continuous function $f:\BR^d\to\BR$ at a point $\vx\in\BR^d$ is defined by
$\partial f(\vx) := {\rm conv}\big\{\vg : \vg = \lim_{\vx_k \rightarrow \vx} \nabla f(\vx_k)\big\}$.
\end{definition}

\begin{definition}[\citet{goldstein1977optimization}]
\label{def:goldsubdiff}
For given $\delta\geq 0$ and a Lipschitz continuous function $f:\BR^d\to\BR$, the Goldstein $\delta$-subdifferential of at point $\vx\in\BR^d$ is defined by
$\partial_{\delta}f(\vx) := {\rm conv}\big(\cup_{\vy \in \BB^d(\vx,\delta)} \partial f(\vy)\big)$,
where the $\partial f(\vy)$ is Clarke subdifferential. 
\end{definition}

We are interested in finding the $(\delta,\epsilon)$-Goldstein stationary point \citep{zhang2020complexity}, which is defined as follows. 

\begin{definition}[\citet{zhang2020complexity}]\label{def:goldstein}
For given Lipschitz continuous function $f:\mathbb{R}^d \rightarrow \mathbb{R}$, $\delta \geq 0$, and~$\vx\in\mathbb{R}^d$, we denote  
$\gradnormdelta{f}{\vx} := \min\{\|\vg\| : \vg \in \partial_{\delta} f(\vx) \}$.
We call the point $\vx\in\BR^d$ as a $(\delta,\epsilon)$-Goldstein stationary point of $f$ if $\gradnormdelta{f}{\vx} \le \epsilon$ holds. 
\end{definition}

\subsection{Randomized Smoothing}

Randomized smoothing is a popular technique in stochastic optimization \citep{duchi2012randomized,lin2022gradient,nesterov2017random,shamir2017optimal,yousefian2012stochastic}.
This paper focuses on the uniform smoothing as follows.

\begin{definition}[\citet{yousefian2012stochastic}] \label{dfn:f-delta}
Given a Lipschitz continuous function $f:\BR^d\to\BR$, we denote its smooth surrogate as
$f_{\delta}(\vx) \triangleq \mathbb{E}_{\vu\sim\un(\BB^d(\vzero,1))}[f(\vx + \delta \vu)]$,
where $\un(\BB^d(\vzero,1))$ is the uniform distribution on the unit Euclidean ball centered at the origin.
\end{definition}

The smooth surrogate $f_\delta$ has the following properties.

\begin{proposition}[{\citet[Proposition 2.2]{lin2022gradient}},  {\citet[Lemma 4]{kornowski2023algorithm}}]\label{prop:ran_smo}
Suppose the function \mbox{$f : \mathbb{R}^d \rightarrow \mathbb{R}$} is $L$-Lipschitz, then its smooth surrogate $f_\delta$ holds: 
(a) $|f_\delta(\cdot) - f(\cdot)| \le \delta L$;
(b) $f_\delta(\cdot)$ is $L$-Lipschitz;
(c) $f_\delta(\cdot)$ is differentiable with $c_0\sqrt{d}L \delta^{-1}$-Lipschitz gradients for some numeric constant $c_0>0$;
(d)~$\partial_{\mu} f_{\delta}(x) \subseteq \partial_{\mu+\delta} f(x)$ for all $\mu\geq 0$. 
\end{proposition}

\section{The Algorithm Design}\label{sec:main}

We propose \underline{d}ecentralized \underline{o}nline-to-nonconvex \underline{c}onversion with \underline{c}lient-data \underline{s}ampling (\DOCS) in Algorithm~\ref{alg:decen}, which incorporates the client-data sampling strategy and the multi-consensus step with Chebyshev acceleration (Algorithm~\ref{Algorithm: fast-gossip})  \citep{ye2023multi,arioli2014chebyshev} into the framework of online-to-nonconvex conversion \citep{cutkosky2023optimal,sahinoglu2024online}. 


Our algorithm design is based on iterating with mini-batch stochastic gradient, where the total sampling budget per iteration is $b$ across all $n$ clients.
Specifically, each iteration allocates the local mini-batch sizes by multinomial distribution, i.e.,
\begin{align} \label{eq:Multinomial}[b_{t,1},\ldots,b_{t,n}]^\top\sim\mathrm{Multinomial}\left(b,\frac{\vone}{n}\right),   
\end{align}
where $b_{t,i}$ is the batch size of the $i$th client at the $t$th iteration.
We then establish the mini-batch stochastic gradient 
\begin{equation}\label{eq:grad_es}
\vg_{t,i}
= \sum_{j=1}^{b_{t,i}}
\nabla F_i(\vw_{t,i}+\mu \vz_{t,i};\vxi_{t,i,j})
\end{equation}
for each client $i\in[n]$,  
where $\vz_{t,i}\sim\un(\BB^d(\vzero,1))$, $\vxi_{t,i,j}\sim\fD_i$, and $\mu\ge 0$.
We then combine the above local gradient estimator and   Chebyshev acceleration (Algorithm~\ref{Algorithm: fast-gossip}) with online-to-nonconvex conversion to establish Algorithm \ref{alg:decen}, where we define 
\begin{align*}
    {\mClip}(\vx,D) = \begin{bmatrix}
        {\vclip}(\vx_1, D) \\
        \vdots \\
        {\vclip}(\vx_n, D)
    \end{bmatrix}\in\BR^{n\times d}
    ~~~\text{with}~~~
    {\vclip}(\vx, D) = \min\left\{1,\frac{D}{\Norm{\vx}}\right\}\vx.
\end{align*}
\begin{remark}
Note that the batch sizes allocation (\ref{eq:Multinomial}) gives 
$\sum_{i=1}^b b_{t,i}=b$ and $\mathbb E[b_{t,i}]=b/n$,
allowing the client with~$b_{t,i}=0$ to skip local stochastic gradient evaluation, i.e., we let $\vg_{t,i}=\vzero$ if $b_{t,i}=0$.    
Since the batch sizes on different clients may be not identical, the $t$th iteration takes the computation rounds of $\max_{i\in[n]} b_{t,i}$, which is not always proportion to the sample complexity $b=\sum_{i=1}^n b_{t,i}$.
\end{remark}

\begin{remark}
The factor $n/b$ in line \ref{line:update-V} of Algorithm \ref{alg:decen} scales the mini-batch stochastic gradients by total sampling budget, which ensures the average of the rows of matrix $n \mG_t/b$ to~be
\begin{align*}
    \hat\vg_t = \frac{1}{b}\sum_{i=1}^n\sum_{j=1}^{b_{t,i}}\nabla F_i(\vw_{t,i}+\mu \vz_{t,i};\vxi_{t,i,j}).
\end{align*}
Since FastGossip (Algorithm~\ref{Algorithm: fast-gossip}) encourages all $\vw_{t,i}$ be close to their average $\bar\vw_{t}=\frac{1}{n}\sum_{i=1}^n \vw_{t,i}$, the vector $\hat\vg_t$ can be regarded as a stochastic gradient evaluation of the global smooth surrogate function~$\frac{1}{n}\sum_{i=1}^n (f_i)_{\mu}(\cdot)$ at the point $\bar\vw_t$ with batch size $b$.
\end{remark}

\begin{algorithm}[t]
\small\caption{Decentralized Online-to-Nonconvex Conversion with Client-Data Sampling (\DOCS)}\label{alg:decen}       
\begin{algorithmic}[1]
    \STATE \textbf{Input:} $\mu\geq 0$,~~$K,T,R,b \in \mathbb{N}$,~~$\eta,\mu,D>0$,~~$\mP\in\BR^{n\times n}$  \\[0.15cm]
    \STATE \textbf{Initialization:} $\mX_{0} = \mV_0 = \vzero$ \\[0.15cm]
    \STATE\textbf{for} $t = 1$ \textbf{to} $KT$ \textbf{do} \\[0.15cm]
    \STATE\quad $\mX_t=\FG(\mX_{t-1}+\mV_t,\mP,R)$ \\[0.15cm]
    \STATE\quad $s_t\sim\un([0,1])$ \\[0.15cm]
    \STATE\quad $\mW_t=\mX_{t-1}+s_t\mV_t$ \\[0.15cm]
    \STATE\quad $\displaystyle{[b_{t,1},\dots,b_{t,n}]^\top\sim {\rm Multinomial}\Big(b,\frac{\vone}{n}\Big)}$ \\[0.15cm] 
    \STATE\quad $\vz_{t,i}\overset{\rm i.i.d}{\sim} \un(\BB^d(\vzero,1))$ for all $i\in[n]$ \\[0.15cm]
    \STATE\quad $\vxi_{t,i,j}\overset{\rm i.i.d}{\sim} \fD_i$ for all $i\in[n]$ and $j\in[b_{t,i}]$ \\[0.15cm]
    \STATE\quad $\displaystyle{\mG_t = \left[\sum_{j=1}^{b_{t,1}}\nabla F_1(\vw_{t,1}+\mu \vz_{t,1};\vxi_{t,1,j});\,\dots\,;\sum_{j=1}^{b_{t,n}}\nabla F_n(\vw_{t,n}+\mu \vz_{t,n};\vxi_{t,n,j})\right]}$ \\[0.1cm]
    \STATE\quad $\displaystyle{\mV_{t+1}={\mClip}\left(\FG\left(\mV_t-\eta\cdot\frac{n\mG_t}{b},\mP,R\right), D\right)}$ \label{line:update-V} \\[0.1cm]
    \STATE \textbf{end for} \\[0.1cm]
    \STATE\textbf{Output:} $\vw_i^{\rm out} \sim \un\big(\big\{\hat{\vw}_i^{(1)}, \dots, \hat{\vw}_i^{(K)}\big\}\big)$ for all $i\in[n]$, where $\hat{\vw}_i^{(k)}=\frac{1}{T}\sum_{t=(k-1)T+1}^{kT}\vw_{t,i}$
\end{algorithmic}
\end{algorithm}

\begin{algorithm}[t]
\small\caption{${\rm FastGossip}(\mZ^{(0)}, \mP, R)$}\label{Algorithm: fast-gossip}
        \begin{algorithmic}[1]
            \STATE \textbf{Initialize:} $\mZ^{(-1)}=\mZ^{(0)}$, $\phi = \dfrac{1-\sqrt{1-(\lambda_2(\mP))^2}}{1+\sqrt{1-(\lambda_2(\mP))^2}}$   \\[0.1cm]
            \STATE \textbf{for} $r = 0$ \textbf{to} $R-1$ \\[0.1cm]
            \STATE ~~~~ $\mZ^{(r+1)} = (1+\phi)\mP \mZ^{(r)} - \phi \mZ^{(r-1)}$ \\[0.1cm]
            \STATE \textbf{end for} \\[0.1cm]
            \STATE \textbf{Output:} $\mZ^{(R)}$
\end{algorithmic}
\end{algorithm}



\section{The Complexity Analysis}\label{sec:complexity}
In this section, we provide a convergence analysis of \DOCS~(Algorithm \ref{alg:decen}). We first consider the case where each local loss function $f_i$ has Lipschitz-continuous gradients and then adopt the analysis to the general nonsmooth nonconvex setting.

\subsection{Analysis for the Smooth Case}





We first introduce the following proposition for the subroutine of multi-consensus steps with Chebyshev acceleration (Algorithm \ref{Algorithm: fast-gossip}) \citep{ye2023multi}.

\begin{proposition}[{\citet[Proposition 1]{ye2023multi}}]\label{prop:chebyshev}
For Algorithm \ref{Algorithm: fast-gossip}, we denote $\bar{\vz}=\frac{1}{n}\vone^\top\mZ^{(0)}$, then it holds that~$\frac{1}{n}\vone^\top\mZ^{(R)}=\bar{\vz}$ and
\begin{align*}
 \|\mZ^{(R)}-\vone\bar{\vz}^\top\|_F^2 \le 14 \left( 1 - \left( 1 - \frac{1}{\sqrt{2}} \right) \sqrt{\gamma} \right)^{2R} \|\mZ^{(0)}-\vone^\top\bar{\vz}\|_F^2.
\end{align*}
\end{proposition}

Based on Proposition \ref{prop:chebyshev}, we can bound the consensus errors of $\mV_t$ and $\mX_t$ by following lemmas.
\begin{lemma}\label{lemma:OL communicate}
Under Assumptions \ref{assumption:gradientoracle} and \ref{assumption:P}, Algorithm \ref{alg:decen} with
\begin{align}\label{eq:R-Delta}
R\geq\left\lceil\frac{1}{(1-1/\sqrt{2})\sqrt{\gamma}}\log\frac{\sqrt{28nD^2+n^2D(G^2+\sigma^2))}}{\epsilon'}\right\rceil
\end{align}
satisfies $\BE[\|\mV_t-\bar{\vv}_{t}\|]\le \epsilon'$ for all $i\in[n]$ and $\epsilon'>0$.
\end{lemma}
\begin{lemma}\label{lem:dec_boundx}
Under the setting of Lemma \ref{lemma:OL communicate}, Algorithm \ref{alg:decen}  holds 
$  \|\mX_t-\bar\vx_t\|\leq \epsilon'$, for all $i \in [n]$ and~$\epsilon' < D$.
\end{lemma}

We then consider the decrease of the objective function value at the mean vector in the smooth case, where each $f_i$ is $H$-smooth.
According to the update rule of $\mX_{t}$ in Algorithm \ref{alg:decen} (line 4) and the doubly stochastic assumption of matrix $\mP$ (Assumption \ref{assumption:P}), we have:
\[\bar{\vx}_t
= \bar{\vx}_{t-1}+\bar\vv_t.\]

Let $\tilde\vnabla_{t}=\int_{0}^{1} \nabla f(\bar{\vx}_{t-1} + s\bar\vv_{t}){\rm d}s
$, 
then the continuity of the function $f$ means
\[f(\bar{\vx}_t) - f(\bar{\vx}_{t-1}) 
= \int_{0}^{1} \langle \nabla f(\bar{\vx}_{t-1} + s\bar\vv_{t}),\bar\vv_{t} \rangle {\rm d}s 
= \langle \tilde\vnabla_{t},\bar\vv_{t} \rangle. \]
Summing the above equation over $t=(k-1)T,(k-1)T+1,\dots,kT$ and taking expectation yields
\begin{align}\label{eq:decompo}
&\BE[f(\bar{\vx}_{kT})- f(\bar{\vx}_{(k-1)T})] = \sum_{t=(k-1)T+1}^{kT}\BE[\langle \tilde\vnabla_{t},\bar\vv_{t} \rangle]\\
=& \underbrace{\sum_{t=(k-1)T+1}^{kT} \BE\left[\left\langle\tilde\vnabla_{t}-\frac{n}{b}\bar\vg_{t},\bar{\vv}_t\right\rangle\right]}_{S_1} 
+ \underbrace{\sum_{t=(k-1)T+1}^{kT} \BE\left[\left\langle \frac{n}{b}\bar\vg_{t},\vu_k\right\rangle\right]}_{S_2}
+ \underbrace{\sum_{t=(k-1)T+1}^{kT} \BE\left[\left\langle \frac{n}{b}\bar\vg_{t},\bar{\vv}_{t}-\vu_k\right\rangle\right]}_{S_3}.
\end{align}
hold for all $\vu_k\in\BR^d$.
We can bound $S_1$, $S_2$ and $S_3$ by the following three lemmas.

\begin{lemma}\label{lem:S1}
Under Assumptions \ref{assumption:Lipschitz}, \ref{assumption:suboptimality}, \ref{assumption:gradientoracle} and \ref{assumption:P}, we further suppose each $f_i$ is $H$-smooth, then Algorithm~\ref{alg:decen} with $\mu=0$ holds
\begin{align*}\label{2nd term}
\sum_{t=(k-1)T+1}^{kT} \BE\left[\left\langle\tilde\vnabla_{t} - \frac{n}{b}\bar\vg_{t}, \bar{\vv}_{t} \right\rangle\right]
\le 2DH\epsilon'T.
\end{align*}
\end{lemma}

\begin{lemma}\label{lem:S2}
Under Assumptions \ref{assumption:Lipschitz}, \ref{assumption:suboptimality}, \ref{assumption:gradientoracle} and \ref{assumption:P}, we let $\vu_k = -D\frac{\sum_{t=(k-1)T+1}^{kT} \finitesum{i}{n} \nabla f_i(\vw_{t,i}) }{\norm{\sum_{t=(k-1)T+1}^{kT} \finitesum{i}{n} \nabla f_i(\vw_{t,i})}}$. Then Algorithm~\ref{alg:decen} with $\mu=0$ holds
{\small \begin{align*}
\BE\left[\left\langle \sum_{t=(k-1)T+1}^{kT}\frac{n}{b}\bar{\vg}_t, \vu_k \right\rangle\right]\le -D\BE\left[\norm{\frac{1}{n}\sum_{t=(k-1)T+1}^{kT} \finitesum{i}{n} \nabla f_i(\vw_{t,i})}\right]+\frac{2\sqrt{2}D(\sigma+2L)\sqrt{T}}{\sqrt{b}}.    
\end{align*}}
\end{lemma}

\begin{lemma} \label{lem:S3}
Under the setting of Lemma \ref{lem:S2}, Algorithm~\ref{alg:decen} with $\mu=0$ holds
    \begin{align*}
        \sum_{t=(k-1)T+1}^{kT} \BE\left[\left\langle \frac{n}{b}\bar\vg_{t},\bar{\vv}_{t}-\vu_k\right\rangle\right]\le D\sqrt{T(L^2+\sigma^2/b)}+\epsilon'T\sqrt{T(L^2+\sigma^2/b)}
    \end{align*}
\end{lemma}

We can also bound the difference between the gradients of the global and local functions as follows.

\begin{lemma}[{\citet[Lemma 2]{sahinoglu2024online}}]
\label{lem:global bound}
Let the functions $\{f_i\}_{i=1}^{n}$ be $H$-smooth and $f(\vx) = \frac{1}{n} \sum_{i=1}^n f_i(\vx)$. 
For given sequence $\{\vw_{t,i}\}_{t,i=1}^{T,n}$, we suppose there exists some $r>0$ such that $\|\vw_{t,i} -\bar{\vw}_t\| \le r$ for all~$i\in[n]$, $k\in[K]$ and $t\in[T]$, then it holds
\begin{align*}
\norm{ \frac{1}{nT} \sum_{t=(k-1)T+1}^{kT} \sum_{i=1}^{n} \nabla f(\mathbf{\vw}_{t,i}) }
\leq
\norm{ \frac{1}{nT} \sum_{t=(k-1)T+1}^{kT} \sum_{i=1}^{n} \nabla f_i(\mathbf{\vw}_{t,i}) }
+ 2rH,
\end{align*}
where $\bar{\vw}_t = \frac{1}{n} \sum_{i=1}^{n} \vw_{t,i}$.
\end{lemma}

Combing all above results, we establish the theoretical guarantee for our method in the smooth case.

\begin{lemma}
\label{lemma:smooth condition}
Under the settings of Lemma \ref{lem:S1}, running Algorithms \ref{alg:decen} with
$\mu=0$,
$T = \fO((L^2+\sigma^2/b)\epsilon^{-2})$,
$K = \fO(\delta^{-1}\epsilon^{-1}\Delta)$,
$D = \fO(\delta/T)$,
$R = \tilde\fO(\gamma^{-1/2})$,
and $\eta=D/\sqrt{T(L^2+\sigma^2/b)}$ can
output $\{\vw^{\rm out}_i\}_{i=1}^n$ such that $\BE\left[\|\nabla f(\vw^{\rm out}_i)\|_\delta\right] \leq \epsilon$ holds for all $i\in[n]$.
\end{lemma}

\subsection{Analysis for the Nonsmooth Case}

For nonsmooth case, we analyze our algorithm by considering its smooth surrogate function. 
Specifically, Proposition~\ref{prop:ran_smo} states the function $f_\delta$ has $H$-Lipschitz continuous gradients with $H=cL\sqrt{d}\delta^{-1}$. 
Additionally, we can connect the Goldstein subdifferential of a Lipschitz continuous function to that of its smooth surrogate as follows.

\begin{lemma}[{\citet[Proposition 2]{sahinoglu2024online}}]
\label{lemma:deltanorm}
From Proposition \ref{prop:ran_smo} and Definition \ref{def:goldstein}, we have
$\|\nabla f(\vx)\|_{\delta} \leq \|\nabla f_{a\delta}(\vx)\|_{(1-a)\delta}$ for all $a \in (0,1)$.
\end{lemma}

Connecting the above lemma with Lemma \ref{lemma:smooth condition}, we can establish our main results for the nonsmooth case as follows.

\begin{theorem}
\label{thm:first-order}
Under Assumptions \ref{assumption:Lipschitz}, \ref{assumption:suboptimality}, \ref{assumption:gradientoracle}, and \ref{assumption:P}, Algorithm \ref{alg:decen} with  
$\mu=\delta/2$, 
$T = \fO((L^2+\sigma^2/b)\epsilon^{-2})$,
$K = \fO(\delta^{-1}\epsilon^{-1}\Delta)$,
$D = \fO(\delta/T)$,
$R = \tilde\fO(\gamma^{-1/2})$,
and $\eta=D/\sqrt{T(L^2+\sigma^2/b)}$
can output $\{\vw^{\rm out}_i\}_{i=1}^n$ such that
$\BE\left[\|\nabla f(\vw^{\rm out}_i)\|_\delta\right] \leq \epsilon$ for all $i\in [n]$.
\end{theorem}

\begin{corollary}\label{coro:1st}
Following the setting of Theorem \ref{thm:first-order}, each client can achieve an $(\delta,\epsilon)$-Goldstein stationary point within the overall sample complexity of $\fO((bL^2+\sigma^2)\Delta\delta^{-1}\epsilon^{-3})$
, the computation rounds of $\tilde\fO((bL^2+\sigma^2)\Delta n^{-1}\delta^{-1}\epsilon^{-3})$
, and the communication rounds of $\tilde\fO((L^2+b^{-1}\sigma^2)\Delta\gamma^{-1/2}\delta^{-1}\epsilon^{-3})$.
\end{corollary}

\begin{remark}\label{remark:batch-size}
The mini-batch size $b$ provides a trade-off between the complexity of computation and  communication. 
Specifically, increasing $b$ reduces the variance-dependent term in the communication complexity from $\sigma^2$ to $\sigma^2/b$, 
In particular, choosing $b=\Theta(\sigma^2/L^2)$ balances the terms $bL^2$ and $\sigma^2$ in the sample complexity and the computation rounds.
Notably, even with $b=1$, our method still achieves better complexity bounds than existing methods~\citep{sahinoglu2024online} in terms of all three measures.
\end{remark}
\begin{remark}
Besides the sharper complexity bounds than ME-DOL \citep{sahinoglu2024online}, the proposed \DOCS~also guarantees every client can achieve a $(\delta,\epsilon)$-Goldstein stationary point in expectation. 
Recall that the theoretical analysis of ME-DOL  only indicates that there exists some point \mbox{$\bar\vw^k=\frac{1}{nT}\sum_{t=1}^T\sum_{i=1}^n\vw^{k,t}_i$} which is an $(\delta,\epsilon)$-Goldstein stationary point, where $\vw_i^{k,t}$ is generated on the $i$th client.
However, achieving such mean vector $\bar\vw^k$ is non-trivial in practice for the decentralized setting. 
In contrast, the point~$\vw_i^{\rm out}$ in our algorithm and theory only depends on its local variables.
\end{remark}

\section{Experiments}

This section empirically compare our \DOCS~with baseline methods ME-DOL \citep{shahrampour2017distributed} by conduction experiments on the following two models:

\textbf{Nonconvex SVM with Capped-$\ell_1$ Penalty:} We first consider the model of nonconvex penalized SVM with capped-$\ell_1$ regularizer \citep{zhang2010analysis} on dataset ``a9a'' and ``cod-rna'', which targets to train the binary classifier $\vx \in \BR^d$ on dataset $\{(\va_i,b_i)\}_{i=1}^m$, where $\va_i \in \BR^d$ and ${b_i} \in \{-1,1\}$ are the feature vector and label for the $i$-th sample. 
We formulate this problem as the following nonsmooth nonconvex problem
\begin{align*}
    \min_{\vx\in\BR^d} f(\vx) \triangleq \frac{1}{m} \sum_{i=1}^m g_i(\vx), 
\end{align*}
where $g_i(\vx)=l(b_i \va_i^\top \vx) +  \nu(\vx)$, $l(z) = \max\{1-z,0\}$, $\nu(\vx)=\lambda \sum_{j=1}^d \min\{|x(j)|,\alpha\}$, and $\lambda,\alpha >0$. 
Here, the notation $x(j)$ means the $j$th coordinate of $\vx$.
We evenly partition $\{g_i\}_{i=1}^m$ into $m$ clients.
and set $\lambda = 10^{-5}/m$ and $\alpha =2$ in our experiments. 

\textbf{Multilayer Perceptron with ReLU activation:}
We then consider the applications of image classification on datasets ``MNIST'' and ``fashion-MNIST'' ($28\times28$ pixels for each image, 10 classes).
Specifically, we consider the two-layer Multilayer Perceptron (MLP) with ReLU activation and a 256-dimensional hidden layer.
Specifically, the local function at the $i$-th client can be written as
\begin{align*}
f_i(\vx) \triangleq \frac{1}{m_i}\sum_{j=1}^{m_i} \ell\big(g(\vx;\va_i^{j}), b_i^{j}\big) + \lambda\|\vx\|_2^2,
\end{align*}
where we organize the parameters of the model as 
$\vx = (\mW_1, \vc_1, \mW_2, \vc_2)$ with $\mW_1 \in \mathbb{R}^{256\times 784}$, $\vc_1 \in \mathbb{R}^{256}$, $\mW_2 \in \mathbb{R}^{10\times 256}$, $\vc_2 \in \mathbb{R}^{10}$ and denote
$g(\vx;\va_i^{j}) = \mW_2 \cdot \text{ReLU}(\mW_1 \va_i^{j} + \vc_1) + \vc_2$ with $\text{ReLU}(x) = \max(0, x)$. 
Additionally, we let
\begin{align*}
\ell(\hat{\vy}, y) = -\sum_{k=0}^{9} \mathbf{1}_{[y=k]}\log \left( \frac{\exp(\hat{y}_{[k]})}{\sum_{j=0}^{9} \exp({\hat y}_{[j]})} \right),
\end{align*}
where $\hat{y}_j$ is the $j$-th coordinate of $\hat{\vy}$.
We also denote $\va_i^j \in \mathbb{R}^{784}$ and $b_i^j \in \{0,1,\dots,9\}$ as the feature (flattened $28\times28$ images) of the $j$th sample on the $i$th client and its corresponding label.


We perform our experiments on a random-walk ring graph with $n=10$ nodes. 
For both ME-DOL and \DOCS, we tune the stepsize $\eta$ and diameter $D$ from $\{0.001, 0.005, 0.01\}$ and $\{0.05, 0.01, 0.005\}$, respectively. 
For our \DOCS, we compare the performance of the algorithm across different mini-batch sizes $b\in\{5,10,20,40,80\}$ and set $R=2$ for the subroutine of Algorithm \ref{Algorithm: fast-gossip}.

\begin{figure}[t]
    \centering
    \begin{tabular}{cccc}
      \includegraphics[scale=0.4]{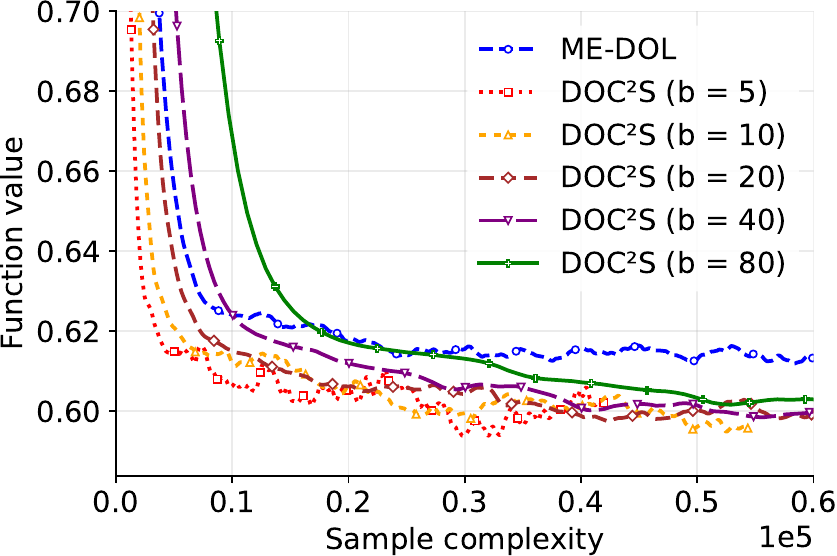}   ~~&~~  
      \includegraphics[scale=0.4]{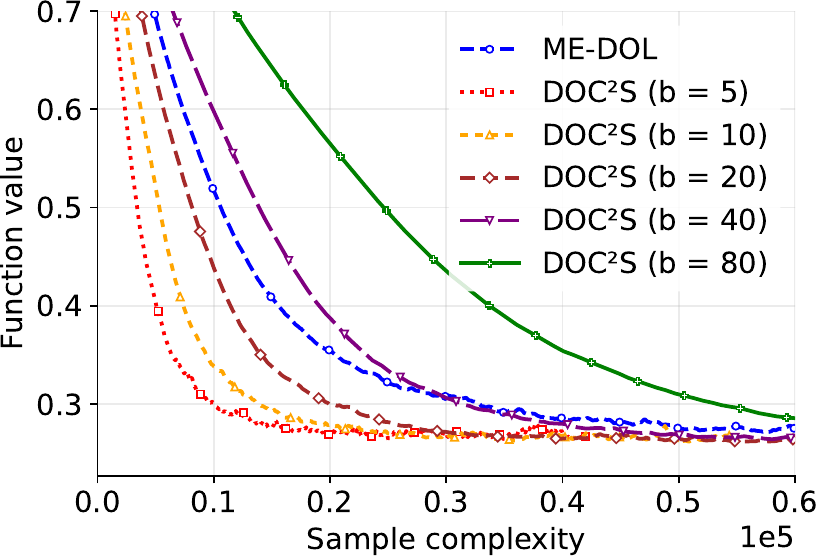}  \\
      (a) a9a & (b) cod-rna \\[0.3cm]
      \includegraphics[scale=0.4]{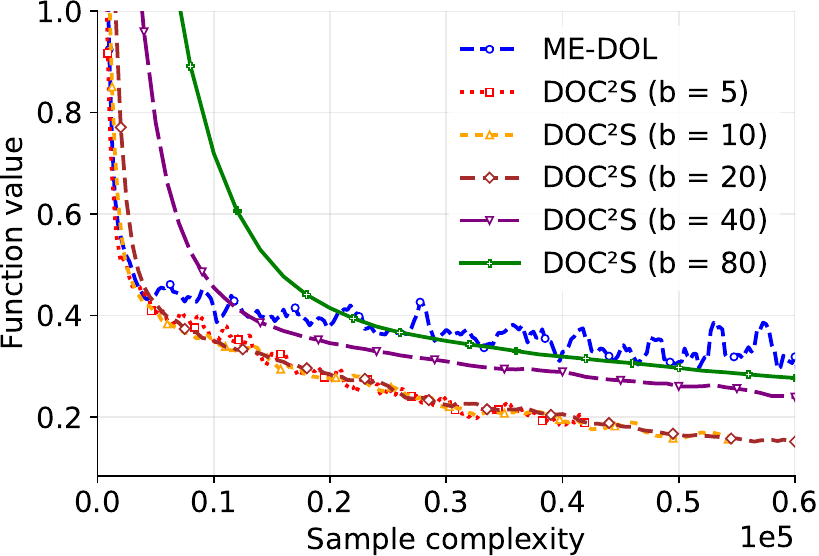} ~~&~~
      \includegraphics[scale=0.4]{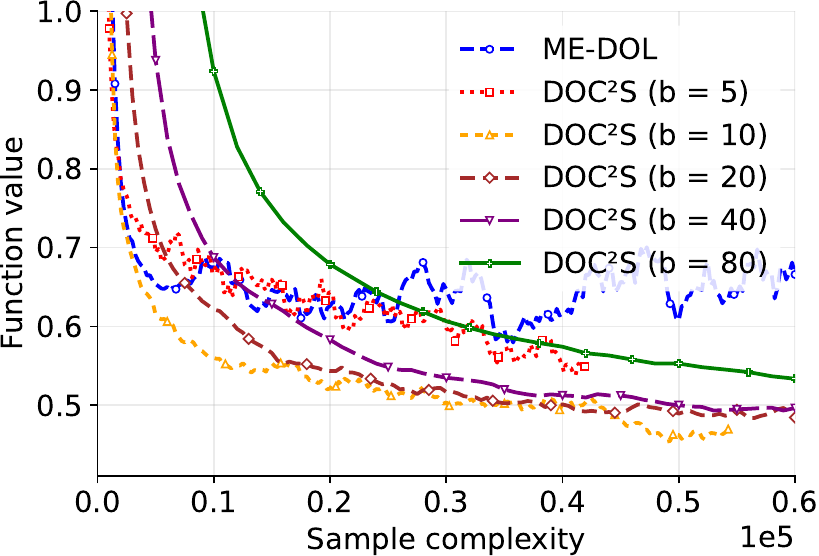} \\
       (c) MNIST & (d) Fashion-MNIST
    \end{tabular}   
    \caption{Comparison on the sample complexity of ME-DOL and DOC$^2$S with different batch sizes on datasets ``a9a'', ``cod-rna'', ``MNIST'', and ``Fashion-MNIST''.}
    \label{fig:sample_complexity}
\end{figure}
\begin{figure}[t]
    \centering
    \begin{tabular}{cccc}
      \includegraphics[scale=0.4]{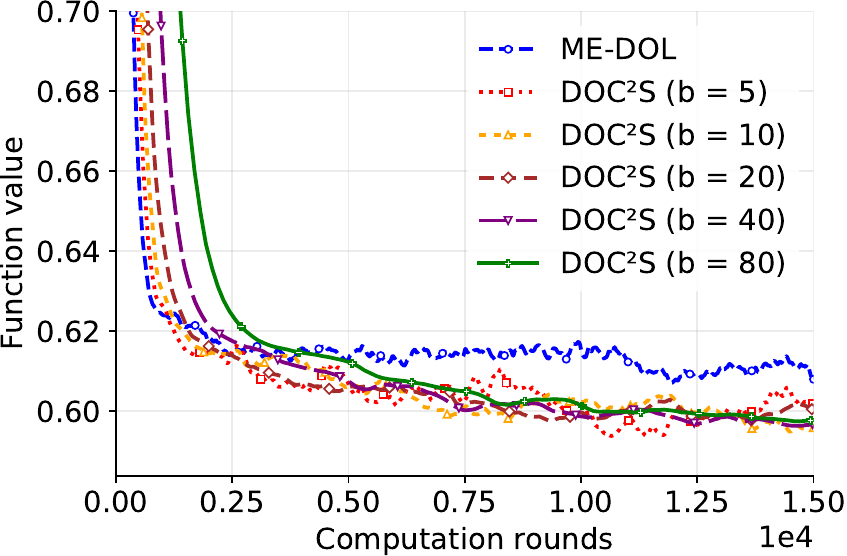}   ~~&~~  
      \includegraphics[scale=0.4]{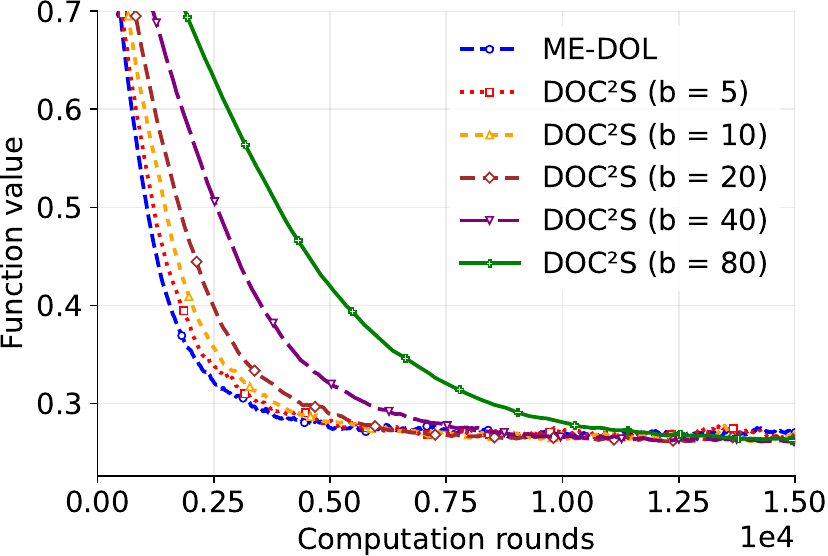} \\ 
      (a) a9a & (b) cod-rna \\[0.3cm]
      \includegraphics[scale=0.4]{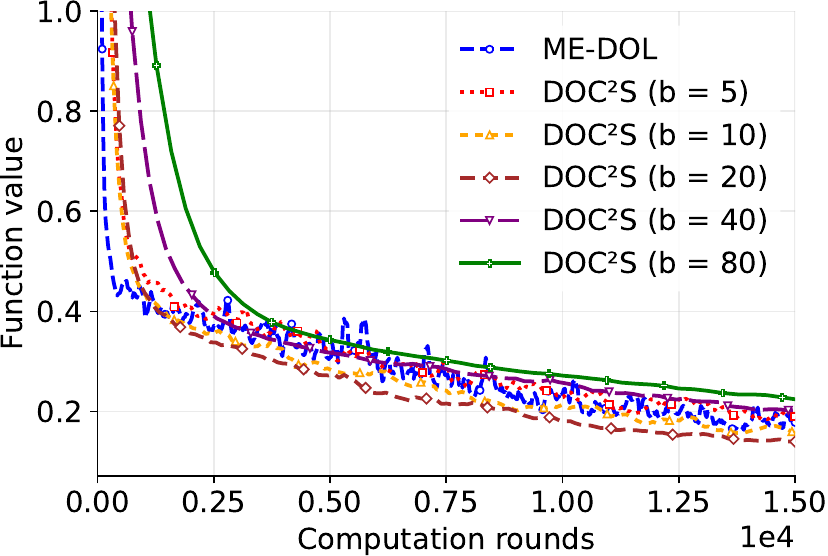} ~~&~~
      \includegraphics[scale=0.4]{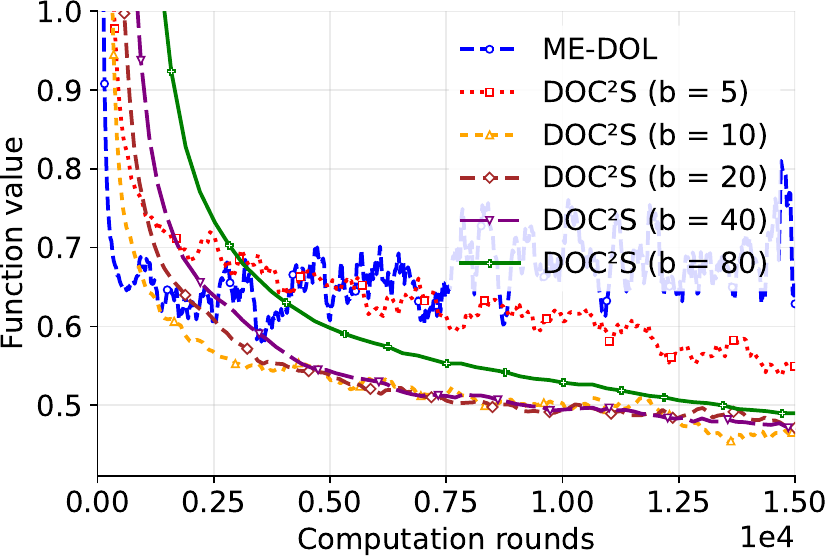} \\
      (c) MNIST & (d) Fashion-MNIST
    \end{tabular}   
    \caption{Comparison on the computation rounds of ME-DOL and DOC$^2$S with different batch sizes on datasets ``a9a'', ``cod-rna'', ``MNIST'', and ``Fashion-MNIST''.}
    \label{fig:computation_rounds}
\end{figure}
\begin{figure}[t]
    \centering
    \begin{tabular}{cccc}
      \includegraphics[scale=0.4]{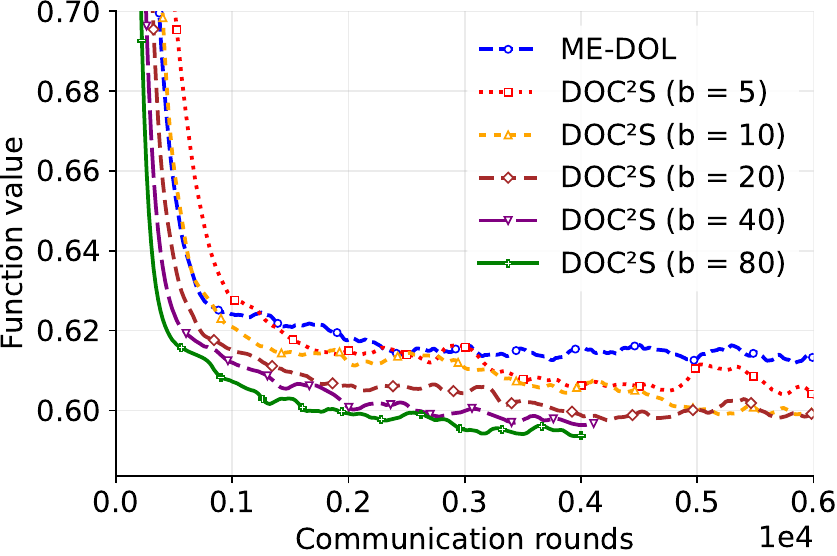}  ~~&~~    
      \includegraphics[scale=0.4]{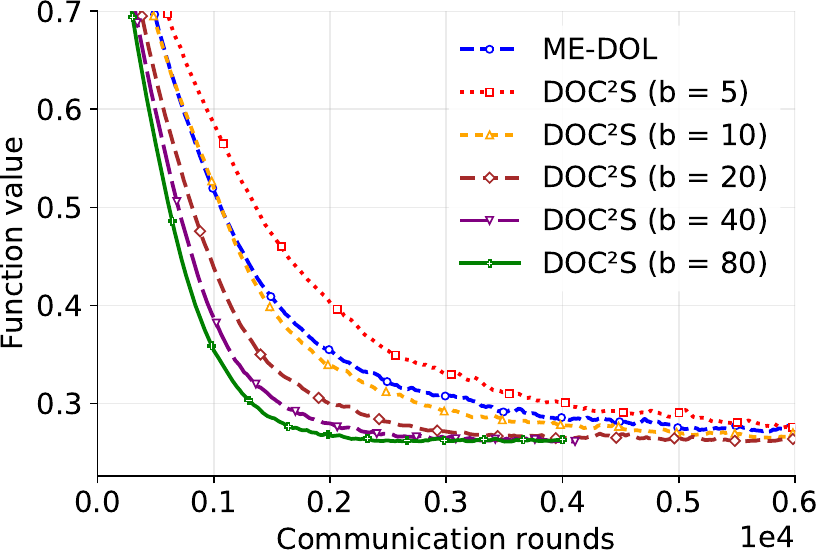} \\
      (a) a9a & (b) cod-rna  \\[0.3cm]
      \includegraphics[scale=0.4]{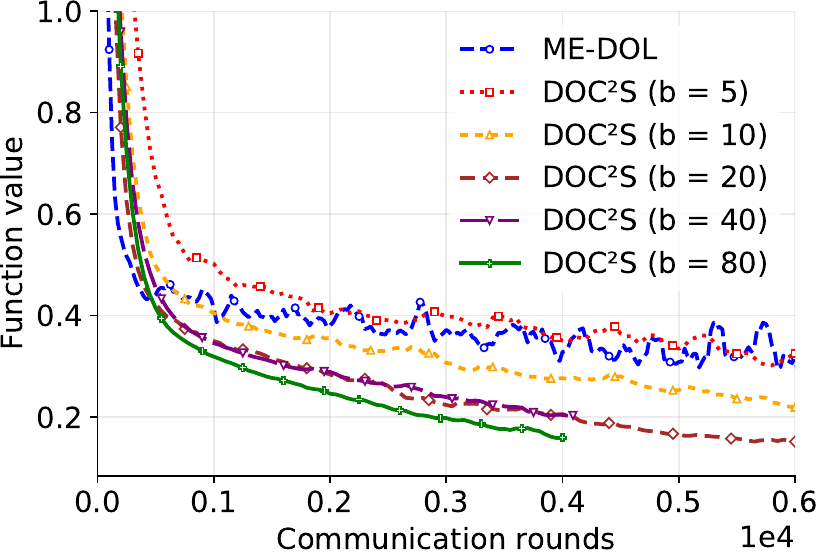} ~~&~~  
      \includegraphics[scale=0.4]{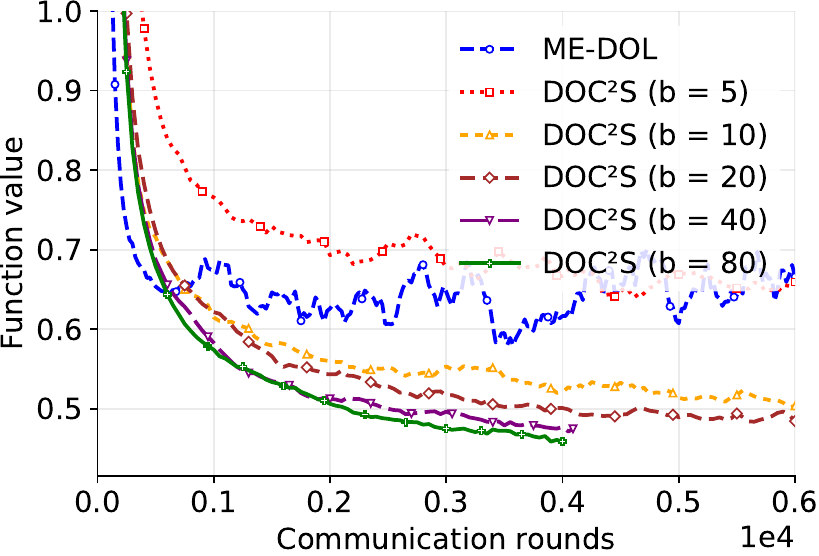} \\
      (c) MNIST & (d) Fashion-MNIST
    \end{tabular}   
    \caption{Comparison on the communication1 rounds of ME-DOL and DOC$^2$S with different batch sizes on datasets ``a9a'', ``cod-rna'', ``MNIST'', and ``Fashion-MNIST''.}
    \label{fig:communication_rounds}
\end{figure}

We present experimental results in terms of sample complexity, computation rounds, and communication rounds in Figures~\ref{fig:sample_complexity}, \ref{fig:computation_rounds},~and~\ref{fig:communication_rounds}, respectively.
We can observe that by taking appropriate mini-batch size, the proposed \DOCS~performs better than baselines with respect to all measures.
For \DOCS, Figure~\ref{fig:communication_rounds} shows that the larger mini-batch size $b$ generally lead to less communication complexity, while the marginal benefit gradually diminishes as $b$ further increases.
This phenomenon supports our theoretical results in Corollary  \ref{coro:1st} and discussion in Remark \ref{remark:batch-size}.



\section{Conclusion}

The paper presents a novel stochastic optimization methods with client-data sampling for decentralized nonsmooth nonconvex problem.
We further provide theoretical analysis to show our method significantly improves the computation and the communication efficiencies.
The advantage of proposed method is also validated by empirical studies. 
In future work, it is possible to extend our ideas to design the efficient zeroth-order methods \citep{lin2022gradient,chen2023faster,lin2023decentralized,kornowski2023algorithm} 
and we can also apply the client-data sampling strategy to address the optmization problem in time varying networks \citep{kovalev2024lower}.

\bibliography{reference.bib}
\bibliographystyle{plainnat}

\newpage
\appendix
\section{Some Useful Lemmas}

We first provide some useful lemmas for our later analysis.

\begin{lemma} \label{lem:sec_m}
    Under Assumptions \ref{assumption:Lipschitz} and \ref{assumption:gradientoracle},  it holds
\begin{align*} 
    \BE[\Norm{\nabla F_i(\vx;\vxi_i)}^2]
\leq 2\BE\left[\Norm{\nabla f_i(\vx)}^2 + \Norm{\nabla f_i(\vx)-\nabla F_i(\vx;\vxi_i)}^2\right]
\leq 2(L^2 + \sigma^2).
\end{align*}

\end{lemma}

\begin{lemma}[{\cite[Lemma 18]{ye2023multi}}] \label{lem:clip}
For all $\mX\in\BR^{n\times d}$ and $D>0$, it holds
\begin{align}
    \Norm{{\mClip}\left(\frac{1}{n}\vone\vone^\top \mX, D\right) - \frac{1}{n}\vone\vone^\top {\mClip}(\mX, D)}_F \leq \Norm{\mX- \frac{1}{n}\vone\vone^\top\mX}_F.    
\end{align}    
\end{lemma}

\begin{lemma}[{\cite[Lemma 4.3]{liu2025nonconvex}}]\label{lemma:heavy-tail}
Given a vector-valued martingale difference sequence $\vv_1,\dots,\vv_T$, we have
\[\BE\left[\left\|\sum_{t=1}^T \vv_t\right\|\right]\le 2\sqrt{2} \BE\left[\left(\sum_{t=1}^T\left\| \vv_t\right\|^p\right)^{\frac{1}{p}}\right], \forall p\in[1,2].\]
\end{lemma}

\begin{lemma}[{\cite[Lemma 6]{shahrampour2017distributed}}]\label{lemma:OL_update}
For given $\vg, \mV \in \mathbb{R}^d$ and $D>0$, the problem
\begin{align}
\min_{\|\vx\|\leq D} \Big\{\langle \vx,\vg \rangle + \frac{1}{2\eta} \|\vx - \mV \|^2\Big\}
\end{align}
has the unique solution 
\begin{align*}
    \vx^*=\min\left\{1,\dfrac{D}{\big\|\mathbf{\Delta} - \eta \vg\big\|}\right\}(\mV-\eta\vg).
\end{align*}
Additionally, we have
\begin{align*}
\langle \vg,\vx^*-\vu \rangle \le \frac{1}{2\eta}\|\mV-\vu\|^2-\frac{1}{2\eta}\|\vu-\vx^*\|^2-\frac{1}{2\eta}\|\mV-\vx^*\|^2,
\end{align*}
for all $\vu\in \mathbb{R}^d$.
\end{lemma}

\section{Upper Bounds for Consensus Errors}

We present the proofs for Lemmas \ref{lemma:OL communicate} and 
\ref{lem:dec_boundx}, which provide upper bounds for consensus errors.

\subsection{Proof of Lemma \ref{lemma:OL communicate}} \label{proof of lemma2}

\begin{proof}
We let 
\begin{align*}
\hat\mV_{t+\frac{1}{2}}=\mV_t-\frac{\eta n\mG_t}{b} 
\qquad\text{and}\qquad 
\mV_{t+\frac{1}{2}}=\FG\left(\hat\mV_{t+\frac{1}{2}},\mP,R\right),    
\end{align*}
then the iteration of Algorithm \ref{alg:decen} and Proposition \ref{prop:chebyshev} implies
\begin{align}\label{eq:delta-consensus-1}
\fnormsq{\mV_{t+\frac{1}{2}}-\vone\bar{\vv}_{t+\frac{1}{2}}^\top} \le 
14c^{2R}\fnormsq{\hat\mV_{t+\frac{1}{2}}-\vone\bar{\vv}_{t+\frac{1}{2}}^\top},
\end{align}
where $c = 1 - ( 1 - {1}/{\sqrt{2})} \sqrt{\gamma}$.
Thus, we have 
\[\BE\left[\fnormsq{\hat\mV_{t+\frac{1}{2}}-\vone\bar{\vv}_{t+\frac{1}{2}}^\top}\right]\le\BE\left[ \fnormsq{\hat\mV_{t+\frac{1}{2}}}\right]=\BE\left[\fnormsq{\mV_t-\frac{\eta n\mG_t}{b}}\right]\le 2n(D^2+\eta n (G^2+\sigma^2)),\]
where the last step comes from $\|\vv_{t,i}\|\le D$ and the fact
\begin{align*}
& \frac{1}{b}\BE[\|\mG_{t}\|^2]
\leq  \frac{1}{b}\sum_{i=1}^{n}\sum_{j=1}^{b_{t,i}}\BE\left[\Norm{
\nabla F_i(\vw_{t,i}+\mu \vz_{t,i};\vxi_{t,i,j})}^2\right] \\
& \leq  \frac{1}{b}\sum_{i=1}^{n}\sum_{j=1}^{b_{t,i}}\BE\left[2\Norm{\nabla f_i(\vw_{t,i}+\mu \vz_{t,i})}^2 +2\Norm{
\nabla F_i(\vw_{t,i}+\mu \vz_{t,i};\vxi_{t,i,j})-\nabla f_i(\vw_{t,i}+\mu \vz_{t,i})}^2  \right] \\
& \leq 2(L^2+\sigma^2)
\end{align*}

By the property of the $\mClip$ operator, we have
\begin{align*}
\fnormsq{\mV_{t+1}-\vone\bar\vv_{t+1}^\top}&=
\fnormsq{\mClip(\mV_{t+\frac{1}{2}},D)-\frac{1}{n}\vone\vone^\top{\mClip}(\mV_{t+\frac{1}{2}},D)}\\
&\le \fnormsq{{\mClip}(\mV_{t+\frac{1}{2}},D)-{\mClip}(\vone\bar\vv_{t+\frac{1}{2}}^\top,D)} + \fnormsq{{\mClip}(\vone\bar\vv_{t+\frac{1}{2}}^\top,D)-\frac{1}{n}\vone\vone^\top{\mClip}(\mV_{t+\frac{1}{2}},D)}\\
&\le 2\fnormsq{\mV_{t+\frac{1}{2}}-\vone\bar\vv_{t+\frac{1}{2}}^\top},
\end{align*}
where the last step is because of Lemma~\ref{lem:clip} and the non-expansiveness of the $\mClip$ operator.

Combining above results, we have
\[\BE\left[\fnorm{\mV_{t+1}-\vone\bar\vv_{t+1}^\top}\right]\le\sqrt{\BE\left[\fnormsq{\mV_{t+1}-\vone\bar\vv_{t+1}^\top}\right]}\le c^{R}\sqrt{28(nD^2+\eta n^2(L^2+\sigma^2))}.\]

Recall that $c = 1 - ( 1 - {1}/{\sqrt{2})} \sqrt{\gamma}$.
Therefore, the setting of $R$ and Proposition \ref{prop:chebyshev} implies

\begin{align}\label{eq:delta-consensus-3}
\BE\left[\fnorm{\mV_{t+1}-\vone\bar\vv_{t+1}^\top}\right]\le \epsilon'  
\end{align}
for all $i\in[n]$.
\end{proof}

\subsection{Proof of Lemma \ref{lem:dec_boundx}}
\begin{proof}
Based on the update rule of $\mX_{t}$ in Algorithm \ref{alg:decen} (line 10), we get
\begin{align}\label{eq:x-recursive}
\begin{split}
& \fnorm{\mX_{t}-\bar{\vx}_{t}} \\
\le& \sqrt{14}c^{R}\fnorm{\mX_{t-1}+\mV_t-(\bar{\vx}_{t-1}+\bar\vv_t)}\\
\le& \sqrt{14}c^{R}\fnorm{\mX_{t-1}-\bar{\vx}_{t-1}} + \sqrt{14}c^{R}\fnorm{\mV_t-\bar\vv_t}\\
\le& \sqrt{14}c^{R}\fnorm{\mX_{t-1}-\bar{\vx}_{t-1}} + \sqrt{14}c^{R}\fnorm{\mV_t}\\
\le& \sqrt{14}c^{R}\fnorm{\mX_{t-1}-\bar{\vx}_{t-1}} + \sqrt{14}c^{R}\cdot \sqrt{n}D,
\end{split}
\end{align}
where the first inequality is based on Proposition \ref{prop:chebyshev}, the second inequality is based on the triangle inequality of Frobenius norm and the last step is based on the fact that $\|\vv_{t,i}\|\le D$ for $i\in[n]$.

We then use induction to prove  
\begin{align*}
    \fnorm{\mX_{t}-\bar{\vx}_{t}}
\leq \epsilon'
\end{align*}
for all $t\geq 1$ and $\epsilon'<D$ as follows.

\paragraph{Induction Base:} For $t=0$, we have 
\begin{align*}
&    \mX_{0}-\bar{\vx}_{0} = 0 \le \epsilon'.
\end{align*}

\paragraph{Induction Step:} We suppose
\begin{align*}
    \fnorm{\mX_{t-1}-\bar{\vx}_{t-1}}
\leq \epsilon'
\end{align*}
holds. Substituting the induction hypothesis  into equation (\ref{eq:x-recursive}) implies
\begin{align*}
    \fnorm{\mX_{t}-\bar{\vx}_{t}}
    \leq \sqrt{14}c^{R}\cdot\epsilon' + \sqrt{14}c^{R}\cdot \sqrt{n}D \leq \epsilon',
\end{align*}
where we take
\begin{align*}
R\geq\left\lceil\frac{1}{(1-1/\sqrt{2})\sqrt{\gamma}}\log\frac{2\sqrt{14n}D}{\epsilon'}\right\rceil.
\end{align*}

\end{proof}

\section{Proofs for the Smooth Case}

We provide the proofs of lemmas for the smooth case in Section \ref{sec:complexity}

\subsection{Proof of Lemma \ref{lem:S1}}
\begin{proof}
Recall that 
\begin{align*}
    \tilde\vnabla_{t}=\int_{0}^{1} \nabla f(\bar{\vx}_{t-1} + s{\bar\vv_{t}})\,{\rm d}s.
\end{align*}
We can split $S_1$ in equation~\ref{eq:decompo} as
\begin{align}\label{decomp_term}
\begin{split}    
&\sum_{t=(k-1)T+1}^{kT} \BE[\langle\tilde\vnabla_{t} - \frac{n}{b}\bar\vg_{t}, \bar{\vv}_{t} \rangle]\\
=& \sum_{t=(k-1)T+1}^{kT} \BE[\langle\bar\vnabla_{t} - \frac{n}{b}\bar\vg_{t}, \bar{\vv}_{t} \rangle]+\sum_{t=(k-1)T+1}^{kT} \BE[\langle\tilde\vnabla_{t} - \bar\vnabla_{t}, \bar{\vv}_{t} \rangle], 
\end{split}
\end{align}
where 
\begin{align*}
\bar{\vnabla}_{t}=\frac{1}{n}\sum_{i=1}^n \vnabla_{t,i}
\qquad\text{and}\qquad
\vnabla_{t,i} =\int_{0}^{1}\nabla f_{i}(\vx_{t-1,i}+s\mV_{t,i})\,{\rm d}s.
\end{align*}
We now consider the upper bound of equation (\ref{decomp_term}). 
Line 6-7 of Algorithm \ref{alg:decen} and Assumption~\ref{assumption:gradientoracle} imply

\[\frac{n}{b}\bar{\vg}_t
=\frac{1}{b}\sum_{j=1}^b \nabla F_{i_{t,j}}(\vw_{t,i_{t,j}};\xi_{t,i_{t,j},k_j})\]
where $i_{t,j} \sim \un(\{1,\dots,n\})$.
We know that 
\begin{align*}
\BE[\frac{n}{b}\bar{\vg}_t]=\BE[\nabla F_{i_{t,j}}(\vw_{t,i_{t,j}};\xi_{t,i_{t,j},k_j})]=\BE[\vnabla_{t,i_{t,j}}]=\BE[\frac{1}{n}\sum_{i=1}^n \vnabla_{t,i}]=\BE[\bar{\vnabla}_{t}],
\end{align*}

Therefore, the first term of equation (\ref{decomp_term}) equal to $0$.

Note that we have
\begin{align}\label{eq:cons-nabla11}
\begin{split}    
&\norm{\tilde\vnabla_{t} - \bar{\vnabla}_{t} } \\
=&\norm{\frac{1}{n}\sum_{i=1}^{n}\int_{0}^{1} \left(\nabla f_i\big(\bar{\vx}_{t-1} + s{\bar\vv_{t}}\big) - \nabla f_{i}\big(\vx_{t-1,i}+s\vv_{t,i}\big)\right)\,{\rm d}s} \\
\leq & \frac{1}{n} \sum_{i=1}^{n} \int_{0}^{1} \norm{\nabla f_i\big(\bar{\vx}_{t-1} + s{\bar\vv_{t}}\big) - \nabla f_{i}\big(\vx_{t-1,i}+s\vv_{t,i}\big)\big)}{\rm d}s\\
\le& \frac{H}{n} \sum_{i=1}^{n} \int_{0}^{1} \norm{\bar{\vx}_{t-1} + s\bar\vv_{t} - \vx_{t-1,i} - s \vv_{t,i}} {\rm d}s \\
\le& \frac{H}{n} \sum_{i=1}^{n} \norm{\bar{\vx}_{t-1} - \vx_{t-1,i}} + \frac{H}{2n} \sum_{i=1}^{n} \norm{\bar\vv_{t} - \vv_{t,i}},
\end{split}
\end{align}
where the second inequality is based on 
the $H$-smoothness of the function $f_i$.

According to Lemma \ref{lem:dec_boundx}, we have
\begin{align}\label{eq:con-y11}
\BE[\|\bar{\vx}_{t-1} - \vx_{t-1,i}\|]\leq \epsilon'.
\end{align}

According to Lemma \ref{lemma:OL communicate}, we have
\begin{align} \label{eq:v_bound}
\BE[\|\bar{\vv}_{t} - \vv_{t,i}\|]\leq \epsilon'
\end{align}

Substituting equations 
(\ref{eq:con-y11}) and (\ref{eq:v_bound}) into equation (\ref{eq:cons-nabla11}), we have
\begin{align}\label{nabla bound 2}
\BE[\|\tilde\vnabla_{t} - \bar{\vnabla}_{t} \|]
\le 2H\epsilon'.
\end{align}

Therefore, the second term in equation (\ref{decomp_term}) holds
\begin{align*}
&\sum_{t=(k-1)T+1}^{kT} \BE[\langle\tilde\vnabla_{t} - \frac{n}{b}\bar\vg_{t}, \bar{\vv}_{t} \rangle]\\
=&\sum_{t=(k-1)T+1}^{kT} \BE[\langle\tilde\vnabla_{t} - \bar{\vnabla}_{t}, \bar{\vv}_{t} \rangle] \\
\le& \sum_{t=(k-1)T+1}^{kT} \BE[\|\tilde\vnabla_{t} - \bar{\vnabla}_{t}\|\|\bar{\vv}_{t}\|]\\
\le& 2DH\epsilon'T,
\end{align*}
where the second step is based on Cauchy--Schwarz inequality, and the last step is based on \eqref{nabla bound 2} and the fact $\|\bar{\vv}_{t}\|\le D$.
\end{proof}

\subsection{Proof of Lemma \ref{lem:S2}}
\begin{proof}
By setting
\begin{align*}
\vu_k = -D\frac{\sum_{t=(k-1)T+1}^{kT} \finitesum{i}{n} \nabla f_i(\vw_{t,i}) }{\norm{\sum_{t=(k-1)T+1}^{kT} \finitesum{i}{n} \nabla f_i(\vw_{t,i})}},
\end{align*}
we have
\begin{align*}\label{eq5}
&\BE\left[\left\langle \sum_{t=(k-1)T+1}^{kT}\frac{n}{b}\bar{\vg}_t, \vu_k \right\rangle\right]\\
=& \BE\left[\left\langle \frac{1}{n}\sum_{t=(k-1)T+1}^{kT}\sum_{i=1}^n \nabla f_i(\vw_{t,i}),\vu_k\right\rangle \right] +\BE\left[ \left\langle \sum_{t=(k-1)T+1}^{kT}\left(\frac{n}{b}\bar{\vg}_t - \frac{1}{n}\finitesum{i}{n} \nabla f_i(\vw_{t,i})\right), \vu_k \right\rangle \right]\\
\le& -D\BE\left[\norm{\frac{1}{n}\sum_{t=(k-1)T+1}^{kT} \finitesum{i}{n} \nabla f_i(\vw_{t,i})}\right] + D\BE\left[\norm{\sum_{t=(k-1)T+1}^{kT} \left(\frac{n}{b}\bar{\vg}_t - \frac{1}{n}\finitesum{i}{n} \nabla f_i(\vw_{t,i})\right) }\right]\\
\le& -D\BE\left[\norm{\frac{1}{n}\sum_{t=(k-1)T+1}^{kT} \finitesum{i}{n} \nabla f_i(\vw_{t,i})}\right]
+ D\BE\left[\norm{\frac{1}{b}\sum_{t=(k-1)T+1}^{kT} \sum_{j=1}^b(\nabla f_{i_{t,j}}(\vw_{t,i_{t,j}}) - \frac{1}{n}\finitesum{i}{n} \nabla f_i(\vw_{t,i})) }\right]\\
&+ D\BE\left[\norm{\frac{1}{b}\sum_{t=(k-1)T+1}^{kT} \sum_{j=1}^b(\nabla F_{i_{t,j}}(\vw_{t,i_{t,j}},\xi_{t,i_{t,j},j}) - \nabla f_{i_{t,j}}(\vw_{t,i_{t,j}})) }\right]
\end{align*}
where the last inequality is based on Cauchy--Schwarz inequality and the fact that $\norm{\vu^k}\leq D$. Since $\BE_{\xi}[\nabla F_{i_{t,j}}(\vw_{t,i_{t,j}})] = \nabla f_{i_{t,j}}(\vw_{t,i_{t,j}})) $ and $\BE_{i_{t,j}}[\nabla f_{i_{t,j}}(\vw_{t,i_{t,j}})]=\frac{1}{n}\finitesum{i}{n} \nabla f_i(\vw_{t,i})$
By Lemma~\ref{lemma:heavy-tail} and Jensen's inequality, we can get
\begin{align*}
&\BE\left[\norm{\sum_{t=(k-1)T+1}^{kT} \sum_{j=1}^b (\nabla F_{i_{t,j}}(\vw_{t,i_{t,j}}) -  \nabla F_{i_{t,j}}(\vw_{t,i_{t,j}})) }\right]\\
\leq& 2\sqrt{2}\mathbb{E}\left[\left(\sum_{t=(k-1)T+1}^{kT}\sum_{j=1}^b\left\|  \nabla F_{i_{t,j}}(\vw_{t,i_{t,j}}) -  \nabla F_{i_{t,j}}(\vw_{t,i_{t,j}})
\right\|^2\right)^{\frac{1}{2}}\right] \\
\leq& 2\sqrt{2}\left(\BE\left[\sum_{t=(k-1)T+1}^{kT}\sum_{j=1}^b\left\| \nabla F_{i_{t,j}}(\vw_{t,i_{t,j}}) -  \nabla F_{i_{t,j}}(\vw_{t,i_{t,j}})
\right\|^2\right]\right)^{\frac{1}{2}} \\
\leq&  2\sqrt{2}(Tb\sigma^2)^{\frac{1}{2}}.
\end{align*}
and
\begin{align*}
&\BE\left[\norm{\sum_{t=(k-1)T+1}^{kT} \sum_{j=1}^b (\nabla f_{i_{t,j}}(\vw_{t,i_{t,j}}) - \frac{1}{n}\finitesum{i}{n} \nabla f_i(\vw_{t,i})) }\right]\\
\leq& 2\sqrt{2}\mathbb{E}\left[\left(\sum_{t=(k-1)T+1}^{kT}\sum_{j=1}^b\left\|  \nabla f_{i_{t,j}}(\vw_{t,i_{t,j}}) - \frac{1}{n}\finitesum{i}{n} \nabla f_i(\vw_{t,i}))
\right\|^2\right)^{\frac{1}{2}}\right] \\
\leq& 2\sqrt{2}\left(\BE\left[\sum_{t=(k-1)T+1}^{kT}\sum_{j=1}^b\left\| \nabla f_{i_{t,j}}(\vw_{t,i_{t,j}}) - \frac{1}{n}\finitesum{i}{n} \nabla f_i(\vw_{t,i}))
\right\|^2\right]\right)^{\frac{1}{2}} \\
\leq&  2\sqrt{2}(Tb(2L)^2)^{\frac{1}{2}}.
\end{align*}
Thus we have
\[\BE\left[\left\langle \sum_{t=(k-1)T+1}^{kT}\frac{n}{b}\bar{\vg}_t, \vu_k \right\rangle\right]\le -D\BE\left[\norm{\frac{1}{n}\sum_{t=(k-1)T+1}^{kT} \finitesum{i}{n} \nabla f_i(\vw_{t,i})}\right]+\frac{2\sqrt{2}D(\sigma+2L)\sqrt{T}}{\sqrt{b}}. \]
\end{proof}

\subsection{Proof of Lemma \ref{lem:S3}} 
\begin{proof}

Let $\tilde\vv_{t+1}=\vclip(\bar\vv_{t}-\frac{n\eta\bar\vg_t}{b},D)$. By Lemma \ref{lemma:OL_update} with $\vg=\frac{n}{b}\bar{\vg}_t$, $\mV=\bar\vv_{t}$, $\vu=\vu_k$ and $\vx^*=\tilde\vv_{t+1}$, we have
\begin{align}\label{regret4}    
\langle \frac{n}{b}\bar\vg_{t},\tilde\vv_{t+1} - \vu_k \rangle 
\leq \frac{1}{2\eta} \big\| \bar\vv_{t} - \vu_k \big\|^2 
- \frac{1}{2\eta} \norm{ \vu_k - \tilde\vv_{t+1} }^2 
- \frac{1}{2\eta} \norm{\bar\vv_{t} - \tilde\vv_{t+1}}^2.
\end{align}
Thus we have
\begin{align*}
\langle \frac{n}{b}\bar\vg_{t},\bar\vv_{t} - \vu_k \rangle 
\leq& \frac{1}{2\eta} \big(\big\| \bar\vv_{t} - \vu_k \big\|^2 
-  \big\| \tilde\vv_{t+1}-\vu_k  \big\|^2 \big)
- \frac{1}{2\eta} \norm{\bar\vv_{t} - \tilde\vv_{t+1}}^2+\langle \frac{n}{b}\bar\vg_{t},\bar\vv_{t} - \tilde\vv_{t+1} \rangle   \\
\le& \frac{1}{2\eta} \big(\big\| \bar\vv_{t} - \vu_k \big\|^2 -  \big\| \bar\vv_{t+1}-\vu_k  \big\|^2 \big)
- \frac{1}{2\eta} \norm{\bar\vv_{t} - \tilde\vv_{t+1}}^2+\langle \frac{n}{b}\bar\vg_{t},\bar\vv_{t} - \tilde\vv_{t+1} \rangle  \\
&+\frac{1}{2\eta}\|\tilde\vv_{t+1}-\bar\vv_{t+1}\|\|\tilde\vv_{t+1}+\bar\vv_{t+1}-2\vu_k\|.
\end{align*}
Note that we have \[\langle \frac{n}{b}\bar\vg_{t},\bar\vv_{t} - \tilde\vv_{t+1} \rangle\le \|\frac{n}{b}\bar\vg_{t}\|\|\bar\vv_{t} - \tilde\vv_{t+1}\|\le \frac{\eta}{2}\|\frac{n}{b}\bar\vg_{t}\|^2+\frac{1}{2\eta} \norm{\bar\vv_{t} - \tilde\vv_{t+1}}^2.\]

By Lemma~\ref{lem:clip}, we also have
\[\BE[\|\tilde\vv_{t+1}-\bar\vv_{t+1}\|]\le \BE\left[\fnorm{\mV_{t+\frac{1}{2}}-\vone\bar\vv^\top_{t+\frac{1}{2}}}\right]\le \epsilon'.\]

Combine all above together, we have 
\begin{align*}
&\BE\left[\sum_{t=(k-1)T+1}^{kT}\langle \frac{n}{b}\bar\vg_{t},\bar\vv_{t} - \vu_k \rangle\right]\\
\le& \frac{\|\bar\vv_0-\vu_k\|^2}{2\eta}+\sum_{t=(k-1)T+1}^{kT}\BE[\frac{\eta}{2}\|\frac{n}{b}\bar\vg_{t}\|^2]+\sum_{t=(k-1)T+1}^{kT}\frac{4D\BE[\|\tilde\vv_{t+1}-\bar\vv_{t+1}\|]}{2\eta}\\
\le& \frac{D^2}{2\eta}+\frac{T\eta(L^2+\sigma^2/b)}{2}+\frac{2TD\epsilon'}{\eta}\\
\le& D\sqrt{T(L^2+\sigma^2/b)}+\epsilon'T\sqrt{T(L^2+\sigma^2/b)},
\end{align*}
where the last step comes from the choice of $\eta=\frac{D}{\sqrt{T(L^2+\sigma^2/b)}}$.
\end{proof}

\subsection{Proof of Lemma \ref{lemma:smooth condition}}\label{proof of lemma smooth} 
\begin{proof}
According to the update rule of $\vx_{t}$ in Algorithm \ref{alg:decen} (line 10), we have:
\[\bar{\vx}_t
= \bar{\vx}_{t-1}+\bar\vv_t.\]
where the last step is based on the doubly stochastic assumption of matrix $\mP$ (Assumption \ref{assumption:P}).

Recall that
\begin{align*}
\tilde\vnabla_{t}=\int_{0}^{1} \nabla f(\bar{\vx}_{t-1} + s\bar\vv_{t}){\rm d}s
\end{align*}
then the continuity of the function $f$ means

\[f(\bar{\vx}_t) - f(\bar{\vx}_{t-1}) 
= \int_{0}^{1} \langle \nabla f(\bar{\vx}_{t-1} + s\bar\vv_{t}),\bar\vv_{t} \rangle {\rm d}s 
= \langle \tilde\vnabla_{t},\bar\vv_{t} \rangle. \]
Summing over $t$ and taking expectation yields
\begin{align*}
\BE[f(\bar{\vx}_{kT})- f(\bar{\vx}_{(k-1)T})] =& \sum_{t=(k-1)T+1}^{kT}\BE[\langle \tilde\vnabla_{t},\bar\vv_{t} \rangle]\\
=& \underbrace{\sum_{t=(k-1)T+1}^{kT} \BE[\langle\tilde\vnabla_{t}-\frac{n}{b}\bar\vg_{t},\bar{\vv}_t\rangle]}_{S_1} + \underbrace{\sum_{t=(k-1)T+1}^{kT} \BE[\langle \frac{n}{b}\bar\vg_{t},\vu_k\rangle]}_{S_2}
+ \underbrace{\sum_{t=(k-1)T+1}^{kT} \BE[\langle \frac{n}{b}\bar\vg_{t},\bar{\vv}_{t}-\vu_k\rangle]}_{S_3}.
\end{align*}
hold for all $k\in[K]$ and $\vu_k\in\BR^d$.

Notice that $S_1$, $S_2$ and $S_3$ can be bounded by Lemma~\ref{lem:S1}, Lemma~\ref{lem:S2} and Lemma~\ref{lem:S3}, respectively.
Thus we have 
\begin{align}\label{eq:sm}
\begin{split}
\BE[f(\bar{\vx}_{kT})- f(\bar{\vx}_{(k-1)T})]\le &2DH\epsilon'T-D\BE\left[\norm{\frac{1}{n}\sum_{t=(k-1)T+1}^{kT} \finitesum{i}{n} \nabla f_i(\vw_{t,i})}\right]\\
&+\frac{2\sqrt{2}D(\sigma+2L)\sqrt{T}}{\sqrt{b}}+D\sqrt{T(L^2+\sigma^2/b)}+\epsilon'T\sqrt{T(L^2+\sigma^2/b)}
\end{split}
\end{align}

Next, we target to bound the distance from $\{\vw_{i}^{t,k}\}$ to $\vw_i^{\rm out}$. For all $i,j\in[n]$, $k\in [K]$, and $t_1, t_2\in[(k-1)T+1,kT]$ such that $t_1>t_2$, we have
\begin{align}\label{w bound}
\norm{ \vw_{t_1,i} - \vw_{t_2,j} } 
\leq \norm{ \vw_{t_1,i} - \bar{\vw}_{t_1} } + \norm{ \bar{\vw}_{t_1} - \bar{\vw}_{t_2} } + \norm{ \bar{\vw}_{t_2} - \vw_{t_2,j} }.
\end{align}
The update rule of $\vw_{t,i}$ implies
\begin{align}\label{eq7}
\begin{split}
& \norm{\vw_{t,i}-\bar{\vw}_{t}} 
=  \norm{\vx_{t-1,i}+s_{t}\vv_{t,i}-\bar{\vx}_{t-1}-\frac{1}{n}\sum_{i=1}^n s_{t}\vv_{t,i}} \\
\le& \norm{\vx_{t-1,i}-\bar{\vx}_{t-1}} + s_t\norm{\vv_{t,i}-\bar\vv}  \\
\le& 2\epsilon'
\end{split}
\end{align}
for all $t\in[T]$ and $\epsilon'\le D/2$, where the last inequality is based on the fact $s_{t} \le 1$ consensus error bounds.
Consequently, we have
\begin{align}\label{eq8}
\begin{split}    
& \norm{ \bar{\vw}_{t_1} - \bar{\vw}_{t_2} }  
\le \sum_{t=t_2}^{t_1-1}\norm{ \bar{\vw}_{t+1} - \bar{\vw}_{t} }  
\le  \sum_{t=1}^T \norm{ \bar{\vw}^{k,t+1} - \bar{\vw}^{k,t} } \\
=& \sum_{t=1}^{T-1} \norm{\bar{\vx}_{t}+\frac{1}{n}\sum_{i=1}^n s_{t+1}\vv_{t+1,i}-\bar{\vx}_{t-1}-\frac{1}{n}\sum_{i=1}^n s_{t}\vv_{t,i}} \\
=& \sum_{t=1}^{T-1} \norm{\bar{\vx}_{t-1}+\vv_{t,i}+\frac{1}{n}\sum_{i=1}^n s_{t+1}\vv_{t+1,i}-\bar{\vx}_{t-1}-\frac{1}{n}\sum_{i=1}^n s_{t}\vv_{t,i}} \\
\le& \sum_{t=1}^T \left(\norm{\vv_{t,i}}+\norm{\frac{1}{n}\sum_{i=1}^n s_{t+1}\vv_{t+1,i}} +\norm{\frac{1}{n}\sum_{i=1}^n s_{t}\vv_{t,i}} \right) \\
\le& 3DT,
\end{split}
\end{align}
where we use the update rule of $\vw_{t,i}$ and 
$\vx_i^{k,t}$ and the setting $s_{t}\in[0,1]$.

Substituting equations (\ref{eq7}) and (\ref{eq8}) into equation (\ref{w bound}), we have
\begin{align}\label{w bound2}
\norm{ \vw_i^{k,t_1} - \vw_j^{k,t_2} }  \le 2\epsilon' + 3DT \le \delta,
\end{align}
where the last inequality is based on taking
\begin{align}\label{D_value}
    D = \frac{\delta}{4T}, \qquad
    T>6,
    \qquad\text{and}\qquad 
    \epsilon'\leq \frac{D}{3T}.
\end{align}

Taking the average on inequality~\ref{eq:sm} over $k=1,\dots,K$ and dividing $DT$, we achieve
\begin{align}\label{eq11}
\begin{split}
&\BE\left[ \frac{1}{K} \sum_{k=1}^{K} \norm{\frac{1}{nT} \sum_{t=(k-1)T+1}^{kT} \sum_{i=1}^{n} \nabla f_i( \mathbf{\vw}_{t,i} ) } \right] \\
\le& \frac{\BE[f(\bar{\vx}^{1,0}) - f(\bar{\vx}^{K,T})]}{DKT} + 2H\epsilon' +\frac{4\sigma}{\sqrt{Tb}}+\frac{8L}{\sqrt{T}}+\epsilon'\sqrt{T(L^2+\sigma^2/b)}/D.
\end{split}
\end{align}

Now we start to show the desired approximate stationary point can be achieved at each client. 
According to Lemma \ref{lem:global bound} with $r=2\epsilon'$, we have
\begin{align}\label{local to global}
\left\|\frac{1}{nT} \sum_{t=1}^{T} \sum_{i=1}^{n} \nabla f(\vw_{i}^{k,t}) \right\| \leq \left\|\frac{1}{nT} \sum_{t=1}^{T} \sum_{i=1}^{n} \nabla f_i(\vw_{i}^{k,t}) \right\| + 2H\epsilon'.
\end{align}
Additionally, we have
\begin{align}\label{local to global1}
    \norm{\hat{\vw}_i^{k}-\vw_j^{k,t}}
    = \norm{\frac{1}{T}\sum_{\tau=1}^T\vw_{i}^{k,\tau}-\vw_j^{k,t}} 
    \le \norm{\frac{1}{T}\sum_{\tau=1}^T\left(\vw_{i}^{k,\tau}-\vw_j^{k,t}\right)}
    \le \delta
\end{align}
for all $i,j\in [n]$, $k \in [K]$, and $t \in [T]$, where the first step is based on the setting of $\hat\vw_{i}^k$ (line 11 of Algorithm~\ref{alg:decen}) and the last step is based on 
\eqref{w bound2}. 
Combing above results, we have
\begin{align*}
&\BE[\norm{\nabla f(\vw_{i}^{\rm out})}_{\delta}] \\
= & \BE\left[\frac{1}{K}\sum_{k=1}^K\norm{\nabla f(\hat\vw_i^k)}_{\delta}\right] \\
\le& \BE \left[ \frac{1}{K} \sum_{k=1}^{K} \norm{\frac{1}{nT} \sum_{t=1}^{T} \sum_{i=1}^{n} \nabla f( \vw_i^{k,t} ) } \right] \\
\le& \frac{\BE[f(\bar{\vx}^{1,0}) - f(\bar{\vx}^{K,T})]}{DKT} + 2H\epsilon' +\frac{4\sigma}{\sqrt{Tb}}+\frac{8L}{\sqrt{T}}+\epsilon'\sqrt{T(L^2+\sigma^2/b)}/D.,
\end{align*}
where the first step is based on the setting of $\vw_i^{\rm out}$ (last line of Algorithm \ref{alg:decen}), 
the second step is based on the definition of $\|\nabla f(\cdot)\|_\delta$ (Definition \ref{def:goldstein}),
and the last step is based on using equations (\ref{eq11}) and~(\ref{local to global}).
By setting $T = \fO((L^2+\sigma^2/b)\epsilon^{-2})$,
$K = \fO(\delta^{-1}\epsilon^{-1}\Delta)$,
$D = \fO(\delta/T)$,
$R = \fO(\gamma^{-1/2}\log(\epsilon'^{-1})$,
and $\eta=D/\sqrt{T(L^2+\sigma^2/b)}$ and $\epsilon'=D/(T(L^2+\sigma^2))+H$,
We can guarantee
\begin{align*}
\BE\left[ \norm{\nabla f(\vw_i^{\rm out})}_\delta\right]\le \epsilon.
\end{align*}

\end{proof}

\section{Proofs for the Nonsmooth Case}

This section follows the proof of Lemma \ref{lemma:smooth condition} to achieve the results in the nonsmooth case.

\subsection{Proof of Theorem \ref{thm:first-order}}\label{proof of thm 1}
\begin{proof}
\label{proof:ns_fo}
Recall that the setting of Theorem \ref{thm:first-order} takes $\mu=\delta/2$, 
then Proposition \ref{prop:ran_smo} means the function $f_{\delta/2}$ is $\delta L/2$-Lipschitz and $2c_0 L\sqrt{d}/\delta$-smooth. 
In the view of minimizing the smooth function $f_{\delta/2}$ by Algorithm \ref{alg:decen},
we can follow the first step in the derivation of \eqref{eq:sm} (in the proof of Lemma \ref{lemma:smooth condition}) by replacing $\delta$, $f$, $\sigma$, and $H$ by $\delta/2$, $f_{\delta/2}$, $\sqrt{L^2+\sigma^2}$, and $2c_0 L\sqrt{d}/\delta$, respectively.
This implies that if we set $\mu=\delta/2$, 
$T = \fO((L^2+\sigma^2/b)\epsilon^{-2})$,
$K = \fO(\delta^{-1}\epsilon^{-1}\Delta)$,
$D = \fO(\delta/T)$,
$R = \tilde\fO(\gamma^{-1/2})$,
$\epsilon'=D/(T(L^2+\sigma^2))+2c_0 L\sqrt{d}/\delta$
and $\eta=D/\sqrt{T(L^2+\sigma^2/b)}$, we can get following results  by using Lemma \ref{lemma:deltanorm}
\begin{align*}
\BE\left[ \norm{\nabla f(\vw_i^{\rm out})}_\delta\right]\le \BE[\|\nabla f_{\delta/2}(\vw_i^{\rm out})\|_{\delta/2}] \le \frac{\epsilon}{3}+\frac{\epsilon}{3}+\frac{\epsilon}{3} = \epsilon
\end{align*}
for all $i \in [n]$. 
Hence,  we can achieve the desired $(\delta,\epsilon)$-Goldstein stationary point on each client with  $T=\fO(\epsilon^{-2})$.
\end{proof}

\subsection{Proof of Corollary \ref{coro:1st}}
\begin{proof}
According to the proof of Theorem \ref{thm:first-order}, we achieve an $(\delta,\epsilon)$-Goldstein stationary point of the objective with $KT = \fO(\delta^{-1}\epsilon^{-3})$. Since $b$ stochastic gradients are sampled according to a multinomial distribution in each round, with probability at least $1-\delta'$, the client receiving the largest number of samples receives at most
\[
\frac{b}{n}+\log\frac{n}{\delta'}
\]
stochastic gradients. Thus the computation rounds should be 
\begin{align*}
\tilde\fO\bigg(\dfrac{(bL^2+\sigma^2)\Delta}{n\delta\epsilon^3}\bigg).
\end{align*}
The overall stochastic first-order oracle complexity is \[\fO\bigg(\dfrac{(bL^2+\sigma^2)\Delta}{\delta\epsilon^3}\bigg).\] 
We perform $R$ communication rounds each time and $R=\tilde{\fO}(\gamma^{-1/2})$ from Lemma \ref{lemma:OL communicate}. Thus the communication rounds is \[\tilde\fO\bigg(\dfrac{(L^2+b^{-1}\sigma^2)\Delta}{\sqrt{\gamma}\delta\epsilon^{3}}\bigg)\].
\end{proof}

\end{document}